\documentclass[oneside, reqno, 11pt, a4paper]{amsart}
\usepackage[australian]{babel}

\AtBeginDocument{
  \DeclareSymbolFont{AMSb}{U}{msb}{m}{n}
  \DeclareSymbolFontAlphabet{\mathbb}{AMSb}
}
\DeclareFontFamily{U}{mathx}{\hyphenchar\font45}
\DeclareFontShape{U}{mathx}{m}{n}{<-> mathx10}{}
\DeclareSymbolFont{mathx}{U}{mathx}{m}{n}
\DeclareMathAccent{\widebar}{0}{mathx}{"73}

\usepackage[dvipsnames]{xcolor}
\usepackage{graphicx}
\usepackage{tikz}
\usepackage{caption}
\usepackage{subcaption}

\usetikzlibrary{shapes.misc,matrix,fit,positioning,arrows.meta,decorations.markings,decorations.shapes,decorations.pathreplacing,calc,shapes.geometric,arrows,shapes,patterns,patterns.meta,math}
\tikzset{Matrix/.style={matrix of nodes, font=\footnotesize,text height=1pt, text depth=0.5pt, text width=8.5pt, align=center, column sep=0pt, row sep=0pt, nodes in empty cells}}
\usepackage{pgf}
\usepgflibrary{fpu}

\graphicspath{{../../}{./}}

\usepackage[a4paper, pdftex, left=2cm, top=2cm, right=2cm, bottom=2cm]{geometry}

\usepackage[final]{pdfpages}
\usepackage{changepage}

\usepackage[foot]{amsaddr}
\numberwithin{equation}{section}

\usepackage{url}
\usepackage{hyperref}
\hypersetup{
    breaklinks=true,
    bookmarksopen=true,
    pdftitle={Geometric Multigrid Preconditioners for conforming Finite Element Methods on Arbitrarily Agglomerated Meshes},
    pdfauthor={Santiago Badia, Jordi Manyer}, 
    pdfsubject={}, 
    colorlinks=true,
    linkcolor=black,
    citecolor=blue,
    filecolor=black,
    urlcolor=blue
}
\newcommand{\eqp}[2]{\eqref{#1}$_{#2}$}
\newcommand{\Eqp}[2]{Equation~\eqref{#1}$_{#2}$}

\newcommand{\sect}[1]{Sect.~\ref{#1}}

\usepackage{csquotes}
\usepackage[backend=biber, defernumbers=true, maxbibnames=5, style=numeric-comp, isbn=false, bibencoding=utf8, safeinputenc, url=false, doi=true, giveninits=true]{biblatex}
\usepackage{mathtools}
\usepackage{mathabx}
\usepackage{stmaryrd}
\usepackage{bbm}

\usepackage{float}
\usepackage{framed}
\usepackage{verbatim}
\usepackage{fancyvrb}
\usepackage{booktabs}
\usepackage{longtable}
\usepackage{colortbl}
\usepackage{multirow}
\usepackage{array}
\usepackage{algorithm}
\usepackage{algpseudocode}
\usepackage{cancel}
\usepackage{accents}
\usepackage{mleftright}

\makeatletter
\algdef{S}[IF]{IfNoThen}[1]{\algorithmicif\ #1}
\makeatother
\usepackage[inline]{enumitem}
\usepackage{siunitx}

\newtheorem{theorem}{Theorem}
\numberwithin{theorem}{section}
\newtheorem{remark}[theorem]{Remark}
\newtheorem{lemma}[theorem]{Lemma}
\newtheorem{corollary}[theorem]{Corollary}

\newtheorem{assumption}[theorem]{Assumption}
\newtheorem{definition}[theorem]{Definition}
\usepackage[breakable]{tcolorbox}
\tcbuselibrary{theorems}
\tcbuselibrary{skins}
\tcbset{
	commonstyle/.style={
		theorem style=plain,
		enhanced jigsaw,
		fonttitle=\bfseries,
		fontupper=\itshape,
		halign=justify,
		separator sign=:,
		description delimiters none,
		description font=\bfseries, 
		terminator sign={.\hspace{0.25em}},
		arc=0mm,outer arc=0mm,
		boxrule=0pt,toprule=0pt,bottomrule=0pt,leftrule=0pt,rightrule=0pt,
		titlerule=0pt,toptitle=0pt,bottomtitle=0pt,top=0pt,
		colback=white,coltitle=black,
		boxsep=0pt, bottom=0pt, left=0pt, 
	}
}
\newtcbtheorem[]{myproblem}{Problem}%
{center, commonstyle, fonttitle=\bfseries}{pb}
\newtcbtheorem[]{mydefinition}{Definition}
{center, commonstyle, fonttitle=\bfseries}{mydef}
\newtcbtheorem[]{myassumption}{Assumption}
{center, commonstyle, fonttitle=\bfseries}{myass}

\usepackage{amsmath, amsfonts, amssymb, amscd, bm}
\mathtoolsset{showonlyrefs}

\newcommand{\identity}{\mathbbm{1}}

\newcommand{\ul}[1]{\underline{#1}}
\newcommand{\ol}[1]{\overline{#1}}
\newcommand{\norm}[2]{\|#2\|_{#1}}
\newcommand{\seminorm}[2]{|#2|_{#1}}
\newcommand{\jump}[1]{\llbracket #1 \rrbracket}
\newcommand{\avg}[1]{\{\!\!\{ #1 \}\!\!\}}

\newcommand{\triplevert}[1]{{\left\vert\kern-0.25ex\left\vert\kern-0.25ex\left\vert #1 \right\vert\kern-0.25ex\right\vert\kern-0.25ex\right\vert}}

\newcommand{\Pk}[1]{\mathbb{P}^{#1}}
\newcommand{\IPk}[1]{\mathcal{P}^{#1}}
\newcommand{\lproj}[2]{\pi_{#1}^{#2}}
\newcommand{\rec}[2]{\mathcal{R}_{#1}}  % potential reconstruction; degree argument kept in the source but not displayed (always q+1), as in the companion paper

\newcommand{\dom}{\Omega}
\newcommand{\T}{\mathcal{T}}
\newcommand{\F}{\mathcal{F}}
\newcommand{\interior}{{\rm in}}
\newcommand{\boundary}{{\rm bd}}

\newcommand{\Tl}[1]{\T_{#1}}
\newcommand{\Fl}[1]{\F_{#1}}

\newcommand{\Th}{\Tl{h}}
\renewcommand{\TH}{\Tl{H}}
\newcommand{\Fh}{\Fl{h}}
\newcommand{\FH}{\Fl{H}}

\newcommand{\Fhi}{\Fh^{\interior}}
\newcommand{\Fhb}{\Fh^{\boundary}}

\newcommand{\FHi}{\FH^{\interior}}  % Interior coarse faces
\newcommand{\FHb}{\FH^{\boundary}}  % Boundary coarse faces

\newcommand{\h}{h}
\renewcommand{\H}{H}    % NOTE: clobbers the double-acute accent \H (avoid in bib names)
\renewcommand{\k}{p}    % Fine polynomial degree  (NOTE: clobbers the ogonek accent \k)
\newcommand{\K}{q}      % Coarse polynomial degree

\newcommand{\HT}{H_T}

\newcommand{\NO}{N_O} % Maximum overlapping subdomains

\newcommand{\Uh}{V_h}         % Fine conforming CG space
\newcommand{\UH}{V_H}         % Coarse auxiliary space
\newcommand{\uh}{u_h}
\newcommand{\vh}{v_h}
\newcommand{\wh}{w_h}
\newcommand{\uH}{u_H}
\newcommand{\vH}{v_H}
\newcommand{\lH}{\lambda_{H}}
\newcommand{\mH}{\mu_{H}}

\newcommand{\Wh}{\mathcal{W}_h}               % Oswald-type averaging operator on the fine mesh
\newcommand{\Ih}{\mathcal{I}_h}               % Element-wise Lagrange interpolant onto the fine broken space
\newcommand{\normDG}[1]{\norm{\mathrm{dg},H}{#1}}  % DG energy norm

\newcommand{\normHHO}[1]{\norm{\ul{\mathrm{hho},H}}{#1}}  % HHO hybrid norm
\newcommand{\normHHOloc}[1]{\seminorm{\ul{\mathrm{hho},T}}{#1}} % local: only a seminorm (kernel = constant pairs), hence | . |

\newcommand{\tr}{\gamma}
\newcommand{\trH}{\tr_{H}}

\newcommand{\isc}[1]{\mathcal{U}_{#1}}

\newcommand{\iscH}{\isc{H}}
\newcommand{\iscHhyb}{\ul{\mathcal{U}}_H}

\newcommand{\UHhyb}{\ul{V}_H}
\newcommand{\MH}{M_H}

\newcommand{\uHhyb}{\ul{u}_H}
\newcommand{\vHhyb}{\ul{v}_H}

\newcommand{\JHhyb}{\ul{J}_H}
\newcommand{\JThyb}{\ul{J}_T}

\newcommand{\uThyb}{\ul{u}_T}
\newcommand{\vThyb}{\ul{v}_T}

\newcommand{\eproj}[2]{\pi_{#1}^{1,#2}}  % Local elliptic projector

\newcommand{\Al}[1]{A_{#1}}
\newcommand{\Ah}{\Al{h}}
\newcommand{\AH}{\Al{H}}
\newcommand{\al}[1]{a_{#1}}
\newcommand{\ah}{a}   % fine form = restriction of a; the fine method is conforming
\newcommand{\aH}{\al{H}}

\newcommand{\IHh}{I_{Hh}}       % Prolongation H → h (abstract)
\newcommand{\IHhdg}{I_{Hh}^{\mathrm{dg}}}    % DG prolongation W_h I_h on broken coarse polynomials
\newcommand{\IHhhho}{I_{Hh}^{\mathrm{hho}}}  % HHO prolongation I_{Hh}^{dg} Theta_H on the skeleton space
\newcommand{\IhH}{I_{Hh}^{*}}   % Restriction (adjoint) h → H

\newcommand{\Dl}[1]{D_{#1}}
\newcommand{\Dh}{\Dl{h}}

\newcommand{\Nhi}{\mathcal{N}_h^{\mathrm{int}}}  % interior Lagrangian nodes
\newcommand{\Ksh}{K^{\sharp}}  % covering simplex of an agglomerate
\newcommand{\BH}{B_H}   % coarse-level solver (approximate inverse of A_H)

\usepackage{acronym}
\acrodef{pde}[PDE]{partial differential equation}
\acrodef{fe}[FE]{finite element}
\acrodef{fem}[FEM]{finite element method}
\acrodefplural{fem}[FEMs]{finite element methods}
\acrodef{dof}[DOF]{degree of freedom}
\acrodefplural{dof}[DOFs]{degrees of freedom}
\acrodef{lhs}[LHS]{left hand side}
\acrodef{rhs}[RHS]{right hand side}

\acrodef{hho}[HHO]{Hybrid High-Order}
\acrodef{hdg}[HDG]{Hybridizable Discontinuous Galerkin}
\acrodef{dg}[DG]{Discontinuous Galerkin}
\acrodef{vem}[VEM]{Virtual Element Method}
\acrodef{wg}[WG]{Weak Galerkin}
\acrodef{dpg}[DPG]{Discontinuous Petrov-Galerkin}

\acrodef{dd}[DD]{Domain Decomposition}
\acrodef{gmg}[GMG]{geometric multigrid}
\acrodef{fgmres}[FGMRES]{flexible generalised minimal residual}
\acrodef{pdf}[PDF]{portable document format}
\acrodef{bddc}[BDDC]{balancing domain decomposition by constraints}

\makeatletter
\def\input@path{{../../}}
\makeatother

\title[]{
Geometric Multigrid Preconditioners for conforming Finite Element Methods on Arbitrarily Agglomerated Meshes
}
\date{\today}

\address{$^\dagger$School of Mathematics\\Monash University\\Clayton\\Victoria 3800\\Australia}
\author[S. Badia]{Santiago Badia$^{\dagger}$}
\email{santiago.badia@monash.edu}
\author[J. Manyer]{Jordi Manyer$^{\dagger}$}
\email{jordi.manyer@monash.edu}

\begin{document}

\begin{abstract}
Geometric multigrid methods for conforming finite element discretisations classically require a nested mesh hierarchy, forcing coarse meshes to be conforming, which is unsuitable for general domains. In order to eliminate this geometric constraint, we introduce and analyse a family of multilevel auxiliary space preconditioners. At the first level, a conforming Lagrangian discretisation on a fine simplicial mesh is preconditioned by a non-conforming polytopal discretisation posed on an arbitrarily agglomerated coarse mesh. At the subsequent levels, the coarse problem is solved by a geometric multigrid method for that polytopal discretisation on further agglomerated meshes, so that the whole hierarchy is generated from the fine mesh by agglomeration alone. We consider two polytopal discretisations, a symmetric interior penalty discontinuous Galerkin method and a hybrid high-order method. The prolongation operator at the first level is an element-wise interpolation onto the fine mesh followed by nodal averaging, preceded in the hybrid high-order case by a local potential reconstruction. For point-Jacobi and vertex-star patch additive Schwarz smoothers, we prove that the condition number of the preconditioned system is bounded by $O((\H/\h)^3)$, where $\H/\h$ is the coarsening ratio, and by $O((\H/\h)^2)$ for the DG variant with penalty proportional to $\H/\h$. Numerical experiments in two and three dimensions, on structured and unstructured meshes with agglomerates generated by automatic clustering algorithms, confirm the robustness of the preconditioners.
\end{abstract}

\maketitle

\section{Introduction}
\label{sec:introduction}

\Ac{gmg} methods are among the most efficient solvers for the large, sparse linear systems arising from the discretisation of \acp{pde} \cite{Hackbusch1985,Bramble1987}. By exploiting a hierarchy of discretisations defined on progressively coarser meshes, they achieve optimal computational complexity, solving the system in a number of operations that scales linearly with the number of unknowns.

For classical \acp{fem}, this efficiency comes at the price of a rigid geometric requirement. Classical \ac{gmg} hierarchies consist of nested simplicial or cuboidal meshes, every level being a conforming mesh of the same element type. Such a hierarchy cannot in general be obtained by coarsening a given fine mesh; the union of a few simplices or cuboids is not a simplex or a cuboid. It has to be generated the other way round, by successive refinement of a coarse mesh. The fine mesh then inherits the geometry of its coarsest ancestor, and the domain boundary is resolved only to the accuracy of the coarsest mesh, which is impractical for general domains. Lifting this restriction requires non-nested hierarchies, which are not automatic either, since they involve unstructured mesh generation at every level, and which in turn demand the construction of intricate inter-level transfer operators \cite{Bramble1991}. Since the coarse and fine spaces are no longer nested, the transfer operators are $L^2$-projections between spaces defined on unrelated meshes. Evaluating them requires locating the fine cells overlapped by each coarse cell and integrating over their intersections, which amounts to building an auxiliary supermesh, a delicate and expensive task, especially in parallel.

Polytopal methods, that is methods supporting meshes whose cells are general polytopes with an arbitrary number of faces, do not have this difficulty. Since most agglomerations of a fine mesh yield an admissible coarser one, the hierarchy can be built from a single fine mesh by successive agglomeration. Multilevel solvers of this kind have been developed for \ac{dg} methods \cite{Antonietti2015Multigrid,Antonietti2017Multigrid}, \ac{vem} \cite{Antonietti2023} and \ac{hho} methods \cite{badia2026geometricmultigridsolvershybrid}.
Polytopal methods are, however, more expensive than their simplicial counterparts, whether because of a larger number of \acp{dof}, as in \ac{dg} methods, or because of a more involved formulation and assembly, as in \ac{vem} and \ac{hho} methods.

In this work we introduce a family of multilevel preconditioners. At the first level, a standard conforming Lagrangian \ac{fe} discretisation on the fine mesh is preconditioned by a non-conforming polytopal discretisation posed on an arbitrarily agglomerated coarse mesh. At the subsequent levels, the coarse problem is solved by an agglomeration-based \ac{gmg} method for that polytopal discretisation, for \ac{dg} the one of \cite{Antonietti2015Multigrid,Antonietti2017Multigrid} and for \ac{hho} the one of \cite{badia2026geometricmultigridsolvershybrid}, so that the whole hierarchy is generated from the fine mesh by agglomeration alone. The convergence of these \ac{gmg} solvers has already been analysed in the cited works, so only the first level needs to be analysed here. This is where the difficulty of transferring between a conforming and a non-conforming space lies. We analyse the two-level method for an arbitrary coarse solver, whose accuracy enters through two constants, and the bounds of the cited works then supply these constants for the \ac{gmg} coarse solver.
The fine mesh can thus be produced by a standard mesher, such as \texttt{Gmsh} \cite{gmsh}, and discretised with a conforming method, which is the cheapest option on the finest and most expensive level, while the coarse mesh can be generated by any agglomeration strategy.
The analysis uses the auxiliary space preconditioning theory \cite{Xu1996,Xu1992}, and covers two polytopal coarse discretisations, a symmetric interior penalty Galerkin (SIPG) \ac{dg} method and a statically condensed \ac{hho} method. The \ac{hho} results extend directly to the closely related \ac{hdg} \cite{Cockburn.Di-Pietro} method, and vice versa. In all cases we prove that the preconditioned system has a condition number bounded independently of the mesh size for a fixed coarsening ratio $\H/\h$, and we corroborate the theory with numerical experiments spanning a range of meshes and agglomeration levels.

Auxiliary space methods have classically been used in the opposite direction to the one proposed here, i.e.\ a \emph{non-conforming} fine discretisation is preconditioned by a \emph{conforming} auxiliary space defined on the same mesh. This is the setting of the seminal works on two-level additive Schwarz preconditioners for non-conforming \acp{fem} \cite{Brenner1996,Oswald1996}, later extended to $C^0$ interior penalty discretisations \cite{Brenner2005}. In these works, the transfer between the non-conforming and conforming spaces is provided by nodal averaging operators. Most solvers for non-conforming hybridised discretisations follow the same direction. Early multigrid methods for \ac{hdg} problems recast the skeleton unknowns into bulk functions, so that a standard conforming piecewise-linear multigrid solver could be used \cite{Cockburn2014MultigridHDG}. The more recent \emph{homogeneous} multigrid methods keep the hybrid discretisation at every level \cite{Di_Pietro-Dong-Kanschat-Rupp-2024-Multigrid}, but conformity is still needed. Their analysis rests on a compatibility condition between the hybrid discretisation and the space of continuous piecewise linear functions, which restricts the method to nested, geometrically conforming simplicial meshes and leaves the agglomerated case open. Such solvers give up the freedom in meshing that motivates hybridised and polytopal methods, and the same holds for the \ac{dg} solvers built on conforming auxiliary spaces mentioned above. Here the fine problem is conforming, and it is the non-conformity of the auxiliary space that allows arbitrary agglomeration.

The closest relatives of our approach are the auxiliary space preconditioners of Kalchev and Vassilevski \cite{Kalchev2020a,Kalchev2020b}, in which a conforming \ac{fe} problem is preconditioned via non-conforming reformulations. Our construction differs in that the auxiliary problem is a standard polytopal discretisation, either SIPG or statically condensed \ac{hho}, equipped, in the \ac{hho} case, with high-order face spaces supported on non-planar agglomerate interfaces \cite{Yemm2022ANA}, and in that we prove condition number bounds that explicitly track the coarsening ratio $\H/\h$. Our method can also be read as a geometric, high-order analogue of smoothed aggregation algebraic multigrid (SA-AMG) \cite{Vanek1996,Brezina2011}, where the tentative prolongator and the algebraic prolongator-smoothing step are replaced by an exact polytopal reconstruction on geometric agglomerates followed by Oswald averaging.

The paper is organised as follows. \sect{sec:preliminaries} fixes the model problem and the conforming discretisation. \sect{sec:auxspace} recalls the auxiliary space framework and reduces the condition number of the two-level preconditioner to three assumptions. \sect{sec:analysis-tools} collects the analysis tools on polytopal elements and \sect{sec:smoother} the additive Schwarz smoothers. \sect{sec:dg} and \sect{sec:hho} verify the three assumptions for the \ac{dg} and \ac{hho} coarse spaces, which gives the bound $O((\H/\h)^3)$, and $O((\H/\h)^2)$ for \ac{dg} with penalty proportional to $\H/\h$. \sect{sec:multilevel} discusses the remaining ingredient, the coarse solver itself, for which agglomeration-based \ac{gmg} methods are available in the literature. \sect{sec:results} reports the numerical experiments and \sect{sec:conclusions} concludes.

\section{Model problem and conforming discretisation}
\label{sec:preliminaries}

Let $\dom \subset \mathbb{R}^d$, $d \in \{2,3\}$, be a bounded Lipschitz polytopal domain with boundary $\partial \dom$. Without loss of generality, we assume $\mathrm{diam}(\dom) = 1$. All hidden constants below are stated under this normalisation. We consider the Poisson problem with homogeneous Dirichlet boundary conditions:
\begin{equation}
  \label{eq:poisson}
  \begin{aligned}
    -\Delta u = f \quad \text{in } \dom, \qquad
    u = 0 \quad \text{on } \partial \dom,
  \end{aligned}
\end{equation}
where $f \in L^2(\dom)$ is a given source term.\footnote{The analysis below can be extended to more general boundary conditions and the following analysis readily applies for a forcing term $f \in H^{-1}(\dom)$, but we will stick to the above setting for simplicity.}  The weak formulation reads: find $u \in H^1_0(\dom)$ such that
\begin{equation}
  \label{eq:weak-poisson}
  a(u,v) \doteq (\nabla u, \nabla v)_{\dom}, \qquad
  \ell(v) \doteq (f, v)_{\dom}, \qquad
  a(u,v) = \ell(v) \quad \forall v \in H^1_0(\dom).
\end{equation}
By the Lax-Milgram lemma, problem \eqref{eq:weak-poisson} admits a unique solution.
Let $\Th$ be a conforming simplicial mesh of $\dom$, consisting of a set of disjoint open simplices $t$ with diameters $h_t = \mathrm{diam}(t)$, and let $\h = \max_{t \in \Th} h_t$ denote the mesh size of the mesh $\Th$. We assume that $\Th$ is shape-regular and quasi-uniform in the usual sense. We restrict the presentation to simplicial meshes and Lagrangian elements for simplicity. All results below extend to tensor-product (quadrilateral/hexahedral) meshes equipped with $\mathbb{Q}^{\k}$ Lagrange elements, for which the nodal and trace arguments used in the analysis hold verbatim.

We denote by $\Fh$ the set of faces of the mesh $\Th$, and by $\Fh(t)$ the set of faces of an element $t \in \Th$. We also introduce the sets of interior and boundary faces as $\Fhi = \{ f \in \Fh : f \subset \dom \}$ and $\Fhb = \{ f \in \Fh : f \subset \partial \dom\}$, respectively.

Given a geometric entity $X$, e.g.\ $t$ or a planar face $f$, let $\Pk{k}(X)$ denote the space of polynomials of total degree at most $k$ on $X$. Similarly, given $\mathcal{X}_h$ a collection of geometrical entities, e.g.\ $\Th$, $\Fh$, we define $\Pk{k}(\mathcal{X}_h) = \bigtimes_{X \in \mathcal{X}_h} \Pk{k}(X)$ the broken polynomial space on each entity.
For any measurable set $X$, $(\cdot,\cdot)_X$ represents the $L^2(X)$ inner product and $\norm{X}{\cdot}$ the associated norm. For a collection $\mathcal{X}_h$, these are extended element-wise, $(u,v)_{\mathcal{X}_h} = \sum_{X \in \mathcal{X}_h} (u,v)_X$. Given $r \geq 0$ and a partition $\mathcal{X}_h$ of $\dom$ (in this work, $\Th$ or $\TH$), we use the \emph{broken} norms and seminorms
\begin{equation}
  \label{eq:broken-norms}
  \norm{H^r(\mathcal{X}_h)}{v}^2 \doteq \sum_{X \in \mathcal{X}_h} \norm{H^r(X)}{v}^2,
  \qquad
  \seminorm{H^r(\mathcal{X}_h)}{v}^2 \doteq \sum_{X \in \mathcal{X}_h} \seminorm{H^r(X)}{v}^2,
\end{equation}
defined for functions that are only piecewise $H^r$-regular over $\mathcal{X}_h$, with the convention $\seminorm{H^0(X)}{v} = \norm{L^2(X)}{v}$. We abbreviate $\seminorm{1,X}{u} \doteq \seminorm{H^1(X)}{u}$ and $\seminorm{1,\mathcal{X}_h}{u} \doteq \seminorm{H^1(\mathcal{X}_h)}{u}$.

Given a polynomial degree $\k \geq 1$, we consider the standard Lagrangian conforming space $\Uh$, that is 
\begin{equation}
  \label{eq:conforming-space}
  \Uh = \{ v_h \in H^1_0(\dom) : v_h|_t \in \Pk{\k}(t) \text{ for all } t \in \Th \}.
\end{equation}
The conforming finite element discretisation of \eqref{eq:weak-poisson} reads: find $\uh \in \Uh$ such that
\begin{equation}
  \label{eq:conforming-discretisation}
  \ah(\uh,\vh) = \ell(\vh) \quad \forall \vh \in \Uh.
\end{equation}
We denote by $\Ah : \Uh \to \Uh$ the associated operator, characterised by $\ah(\uh, \vh) = (\Ah\uh, \vh)_{\dom}$ for all $\uh, \vh \in \Uh$.

\section{Auxiliary space preconditioners: abstract framework}
\label{sec:auxspace}

For a given mesh $\Th$, we define a coarse partition $\TH$, given by non-overlapping agglomeration of elements of $\Th$. That is, $\TH = \{ T \}$ is such that each $T \in \TH$ is a union of elements from $\Th$. Here $H_{T} = \mathrm{diam}(T)$ denotes the diameter of each agglomerated element and $\H = \max_{T \in \TH} H_{T}$ the mesh size of the partition.

Throughout this work, we make the following regularity assumptions on the coarse partition.

\begin{assumption}[Coarse mesh regularity]
  \label{ass:coarse-mesh-regularity}
  There exists a constant $c_R > 0$, independent of $\h$ and $\H$, such that the coarse partition $\TH$ satisfies:
  \begin{enumerate}[label=(\roman*)]
    \item \emph{Quasi-uniformity:} $c_R \H \leq H_T \leq \H$ for all $T \in \TH$;
    \item \emph{Star-shapedness:} every $T \in \TH$ is star-shaped with respect to a ball of radius $\rho_T \geq c_R H_T$;
    \item \emph{Simplicial covering:} there exists a covering $\TH^\sharp = \{\Ksh\}$ of $\TH$ by shape-regular $d$-simplices such that, for every $T \in \TH$, there is $\Ksh \in \TH^\sharp$ with $T \subset \Ksh$ and $\mathrm{diam}(\Ksh) \leq c_R^{-1} H_T$.
  \end{enumerate}
\end{assumption}

Conditions (ii)--(iii) are the standard shape-regularity requirements for polytopal meshes \cite{Cangiani2017hpVersion}. Condition (ii) guarantees that the trace and Poincar\'e inequalities of \sect{sec:analysis-tools} hold on the agglomerated elements with constants depending only on $c_R$, and that each agglomerate is connected. Conditions (ii) and (iii) together yield the inverse inequalities for polynomials on polytopal elements required by the coarse bilinear forms introduced in \sect{sec:dg} and \sect{sec:hho}. The inverse inequality holds on the shape-regular simplex $\Ksh$ supplied by condition (iii) and is transferred to $T$ using $T \subset \Ksh$, while condition (ii) is what allows the resulting $L^2(\Ksh)$ norm to be brought back to $T$.

\begin{remark}[Mesh assumptions in practice]
  The agglomerates generated by the automatic clustering strategies used in \sect{sec:results} are not guaranteed to satisfy Assumption~\ref{ass:coarse-mesh-regularity} uniformly. Star-shapedness, for instance, may fail for complicated agglomerate shapes. The numerical experiments therein indicate that the proposed preconditioners remain effective well beyond the setting covered by the analysis.
\end{remark}

We denote by $\FH = \FHi \cup \FHb$ the set of \emph{interfaces} of the coarse partition $\TH$. An \emph{interior} interface $F \in \FHi$ is the union of the fine faces contained in $\ol{T_1} \cap \ol{T_2}$, for two distinct elements $T_1, T_2 \in \TH$, whenever that union has positive $(d-1)$-dimensional measure, and a \emph{boundary} interface $F \in \FHb$ is the union of the fine faces contained in $\ol{T} \cap \partial\dom$ for some $T \in \TH$. Interfaces are piecewise-linear curves for $d=2$ and polygonal surfaces for $d=3$. We will still refer to these interfaces as \emph{faces} for brevity. For each cell $T$, we use $\partial T$ to denote its boundary, and let $\FH(T) = \{ F \in \FH : F \subset \partial T \}$ denote the set of faces of $T$.
For each face $F \in \FH$, we define $H_F$ as the mean of the diameters of the (at most two) coarse elements sharing $F$:
\[
  H_F \doteq \frac{1}{|\mathcal{T}_H(F)|}\sum_{T \in \mathcal{T}_H(F)} H_T, \qquad
  \mathcal{T}_H(F) \doteq \bigl\{T \in \TH : F \subset \partial T\bigr\}.
\]
In particular, $H_F \simeq \H$ for all $F \in \FH$.

Let $\UH$ be an auxiliary finite-dimensional space defined on the coarse partition $\TH$ and $a_H : \UH \times \UH \to \mathbb{R}$ a symmetric positive definite bilinear form defined on $\UH$. Let $\IHh : \UH \to \Uh$ be a prolongation operator from the auxiliary space $\UH$ to the original space $\Uh$. We also introduce a symmetric positive definite smoother $\Dh : \Uh \to \Uh$ and a symmetric positive definite operator $\BH : \UH \to \UH$, which plays the role of an approximate inverse of $\AH$, that is, of a solver for the coarse problem. We denote by $\AH : \UH \to \UH$ the operator associated with $\aH$, characterised by $(\AH \uH, \vH) = \aH(\uH, \vH)$ for all $\uH, \vH \in \UH$, and by $\IHh^* : \Uh \to \UH$ the adjoint of the prolongation operator, characterised by $(\IHh \vH, \vh)_{\dom} = (\vH, \IHh^* \vh)$ for all $\vH \in \UH$ and $\vh \in \Uh$. Here $(\cdot,\cdot)$ denotes an inner product on $\UH$. The coarse correction $\IHh \AH^{-1} \IhH$, and hence the preconditioner itself, is the same for every choice of this inner product, so the choice affects neither the method nor its implementation. Given a symmetric positive definite bilinear form $b$ (resp. operator $B$) on a space $V$, we denote the associated energy norm by $\norm{b}{v} \doteq b(v,v)^{1/2}$ (resp. $\norm{B}{v} \doteq (B v, v)^{1/2}$) for $v \in V$. We will use the energy norms $\norm{\ah}{\cdot}$, $\norm{a_H}{\cdot}$ and $\norm{\Dh}{\cdot}$ below.
Our auxiliary space preconditioner is then defined as
\begin{equation}
  \label{eq:auxspace-preconditioner}
  P_{\mathrm{add}}^{-1} \doteq \IHh \BH \IhH + \Dh^{-1}.
\end{equation}
\begin{remark}
  \label{rem:multiplicative}
  A multiplicative variant of this preconditioner is defined in \sect{sec:multiplicative} and used in the experiments of \sect{sec:results}. 
  We restrict the analysis to $P_{\mathrm{add}}$. Multiplicative and hybrid variants can be treated with the same ingredients, following \cite{Xu1992}, and we do not pursue this here.
\end{remark}

\begin{assumption}[Stability of the prolongation]
  \label{ass:prolongation-stability}
  There exists a constant $c_I$ such that for all $\vH \in \UH$ it holds that
  \[
    \norm{\ah}{\IHh \vH}^2 \leq c_I^2 \norm{a_H}{\vH}^2
  \]
\end{assumption}

\begin{assumption}[Smoother stability]
  \label{ass:smoother-stability}
  There exists a constant $c_D$ such that for all $\vh \in \Uh$ it holds that
  \[
    \norm{\ah}{\vh}^2 \leq c_D^2 \norm{\Dh}{\vh}^2
  \]
\end{assumption}

\begin{assumption}[Stable decomposition]
  \label{ass:stable-decomposition}
  There exist constants $c_{0,H}$ and $c_{0,D}$ such that for all $\vh \in \Uh$ there exists $\vH \in \UH$ with
  \[
    \norm{a_H}{\vH}^2 \leq c_{0,H}^2 \norm{\ah}{\vh}^2
    \qquad\text{and}\qquad
    \norm{\Dh}{\vh - \IHh \vH}^2 \leq c_{0,D}^2 \norm{\ah}{\vh}^2 .
  \]
\end{assumption}

\begin{assumption}[Accuracy of the coarse solver]
  \label{ass:coarse-solver}
  There exist constants $0 < \lambda_0 \leq \lambda_1$ such that for all $\mH \in \UH$ it holds that
  \[
    \lambda_0 (\AH^{-1}\mH, \mH) \leq (\BH \mH, \mH) \leq \lambda_1 (\AH^{-1}\mH, \mH),
  \]
  that is, $\lambda_0$ and $\lambda_1$ bound the extreme eigenvalues of $\BH \AH$ from below and
  above, so that $\lambda_1/\lambda_0$ bounds the condition number $\kappa(\BH\AH)$ of the
  preconditioned coarse problem. An exact coarse solve, $\BH = \AH^{-1}$, corresponds to
  $\lambda_0 = \lambda_1 = 1$.
\end{assumption}

In the following, given symmetric operators $X$ and $Y$, we use the compact notation $X \preceq Y$ (equivalently, $Y \succeq X$) to indicate that $Y - X$ is positive semi-definite. In this notation, Assumption \ref{ass:coarse-solver} reads $\lambda_0 \AH^{-1} \preceq \BH \preceq \lambda_1 \AH^{-1}$. We recall that inversion reverses this order on symmetric positive definite operators.

Then the condition number of the preconditioned system $P_{\mathrm{add}}^{-1} \Ah$ can be bounded as follows:

\begin{lemma}
  \label{lem:auxspace}
  Assume that Assumptions \ref{ass:prolongation-stability}, \ref{ass:smoother-stability}, \ref{ass:coarse-solver} and \ref{ass:stable-decomposition} hold. Then
  \[
    \kappa(P_{\mathrm{add}}^{-1} \Ah) \leq \bigl(c_D^2 + \lambda_1 c_I^2\bigr) \Bigl( c_{0,D}^2 + \frac{c_{0,H}^2}{\lambda_0} \Bigr)
  \]
\end{lemma}

\noindent This result is a variant of \cite[Theorem 2.1]{Xu1996}. It differs from it in the assumptions on the smoother and on the decomposition. The analysis in \cite{Xu1996} requires a two-sided Richardson-type equivalence for the smoother and measures the stable decomposition in a scaled $L^2$ norm, whereas Assumption \ref{ass:smoother-stability} is one-sided for a general symmetric positive definite $\Dh$ and Assumption \ref{ass:stable-decomposition} is stated in the $\Dh$-norm. Since our assumptions are weaker, the result does not follow from the cited theorem, and we must prove it.
\begin{proof}
  Note first that $P_{\mathrm{add}}^{-1}$ is symmetric with respect to $(\cdot,\cdot)$, so $P_{\mathrm{add}}^{-1}\Ah$ is self-adjoint with respect to $\ah(\cdot,\cdot)$, and it is positive definite, since $\ah(P_{\mathrm{add}}^{-1}\Ah\vh, \vh) \geq (\Dh^{-1}\Ah\vh, \Ah\vh) > 0$ for $\vh \neq 0$. Its extreme eigenvalues are therefore characterised by the Rayleigh quotients below.
  Write $g_v \doteq \IhH \Ah \vh$, so that
  \[
    \ah(P_{\mathrm{add}}^{-1}\Ah\vh, \vh) = (\BH g_v, g_v) + (\Dh^{-1}\Ah\vh, \Ah\vh).
  \]
  For the upper bound, the dual-norm characterisation of the $\Dh^{-1}$-norm as the dual of the $\Dh$-norm  followed by Assumption \ref{ass:smoother-stability}, gives
  \[
    (\Dh^{-1}\Ah\vh, \Ah\vh) = \sup_{\wh \neq 0} \frac{(\Ah\vh, \wh)^2}{(\Dh\wh, \wh)} = \sup_{\wh \neq 0} \frac{\ah(\vh, \wh)^2}{\norm{\Dh}{\wh}^2} \leq c_D^2 \, \norm{\ah}{\vh}^2,
  \]
  while Assumptions \ref{ass:coarse-solver} and \ref{ass:prolongation-stability} give
  \[
    (\BH g_v, g_v) \leq \lambda_1 (\AH^{-1} g_v, g_v) = \lambda_1 \sup_{\mu_H \neq 0} \frac{(g_v, \mu_H)^2}{\norm{a_H}{\mu_H}^2} = \lambda_1 \sup_{\mu_H \neq 0} \frac{\ah(\vh, \IHh\mu_H)^2}{\norm{a_H}{\mu_H}^2} \leq \lambda_1 c_I^2 \, \norm{\ah}{\vh}^2,
  \]
  so that $\lambda_{\max}(P_{\mathrm{add}}^{-1}\Ah) \leq c_D^2 + \lambda_1 c_I^2$.
  For the lower bound, let $\vh = \IHh\vH + \wh$ be the decomposition of Assumption \ref{ass:stable-decomposition}, with $\norm{a_H}{\vH}^2 \leq c_{0,H}^2 \norm{\ah}{\vh}^2$ and $\norm{\Dh}{\wh}^2 \leq c_{0,D}^2 \norm{\ah}{\vh}^2$. Then
  \[
    \norm{\ah}{\vh}^2 = \ah(\vh, \IHh\vH) + \ah(\vh, \wh) = (g_v, \vH) + (\Ah\vh, \wh),
  \]
  and the Cauchy--Schwarz inequality, first in the $\BH$- and $\Dh$-inner products and then in $\mathbb{R}^2$, yields
  \[
    \norm{\ah}{\vh}^2 \leq \bigl[ (\BH g_v, g_v) + (\Dh^{-1}\Ah\vh, \Ah\vh) \bigr]^{1/2} \bigl[ (\BH^{-1}\vH, \vH) + (\Dh\wh, \wh) \bigr]^{1/2}.
  \]
  The first bracket is $\ah(P_{\mathrm{add}}^{-1}\Ah\vh, \vh)$. For the second, Assumption \ref{ass:coarse-solver} gives $\BH \succeq \lambda_0 \AH^{-1}$, hence $\BH^{-1} \preceq \lambda_0^{-1}\AH$ by order reversal under inversion, and $(\BH^{-1}\vH, \vH) \leq \lambda_0^{-1} \norm{a_H}{\vH}^2 \leq \lambda_0^{-1} c_{0,H}^2 \norm{\ah}{\vh}^2$, so that
  \[
    \norm{\ah}{\vh}^4 \leq \ah(P_{\mathrm{add}}^{-1}\Ah\vh, \vh) \Bigl( \frac{c_{0,H}^2}{\lambda_0} + c_{0,D}^2 \Bigr) \norm{\ah}{\vh}^2,
  \]
  i.e.\ $\lambda_{\min}(P_{\mathrm{add}}^{-1}\Ah) \geq \bigl( c_{0,D}^2 + c_{0,H}^2/\lambda_0 \bigr)^{-1}$. Since $P_{\mathrm{add}}^{-1}\Ah$ is self-adjoint and positive definite with respect to $\ah(\cdot,\cdot)$, the two bounds give the claim.
\end{proof}

We now provide a more concrete version of the previous lemma, which will be used in the analysis of the preconditioners proposed in this work.

\begin{lemma}
  \label{lem:auxspace-concrete}
  Assume that Assumption \ref{ass:coarse-solver} holds, and that there exist constants $c_D > 0$ and $c_{D,0}, c_{D,1} \geq 0$ such that the following properties hold:
  \begin{itemize}
    \item The smoother $\Dh$ fulfils that for all $\vh \in \Uh$
    \begin{equation}
      \label{eq:auxspace-abstract-smoother}
      (\Dh\vh, \vh)_{\dom} \leq c_{D,0} \norm{L^2(\dom)}{\vh}^2 + c_{D,1} \norm{H^1(\dom)}{\vh}^2 \ , \quad \norm{\ah}{\vh}^2 \leq c_D^2 (\Dh\vh, \vh)_{\dom}
    \end{equation}
    \item The prolongation operator $\IHh$ fulfils that for all $\vH \in \UH$
    \begin{equation}
      \label{eq:auxspace-abstract-prolongation}
      \norm{\ah}{\IHh \vH}^2 \lesssim \frac{\H}{\h} \norm{a_H}{\vH}^2
    \end{equation}
    \item For all $\vh \in \Uh$, there exists $\vH \in \UH$ such that
    \begin{equation}
      \label{eq:auxspace-abstract-decomposition}
      \norm{a_H}{\vH}^2 \lesssim \norm{\ah}{\vh}^2 \ , \quad \norm{H^1(\dom)}{\vh - \IHh \vH}^2 \lesssim \frac{\H}{\h} \norm{\ah}{\vh}^2 \ , \quad \norm{L^2(\dom)}{\vh - \IHh \vH}^2 \lesssim \H^2 \norm{\ah}{\vh}^2
    \end{equation}
  \end{itemize}
  Then the condition number of the preconditioned system satisfies
  \begin{equation}
    \kappa(P_{\mathrm{add}}^{-1} \Ah) \lesssim \left(c_D^2 + \lambda_1 \frac{\H}{\h}\right)\left(c_{D,1} \frac{\H}{\h} + c_{D,0} \H^2 + \frac{1}{\lambda_0}\right).
  \end{equation}
\end{lemma}
\begin{proof}
  Equation~\eqref{eq:auxspace-abstract-prolongation} gives Assumption \ref{ass:prolongation-stability} with $c_I^2 = \frac{\H}{\h}$.
  Equation~\eqp{eq:auxspace-abstract-smoother}{2} is Assumption \ref{ass:smoother-stability} with the constant $c_D$.
  Finally, we combine \eqp{eq:auxspace-abstract-smoother}{1} and \eqref{eq:auxspace-abstract-decomposition}. For all $\vh \in \Uh$, take $\vH \in \UH$ as in \eqref{eq:auxspace-abstract-decomposition} and set $\wh = \vh - \IHh \vH$. Then \eqp{eq:auxspace-abstract-decomposition}{1} bounds the first term of Assumption \ref{ass:stable-decomposition}, while \eqp{eq:auxspace-abstract-decomposition}{2,3} together with \eqp{eq:auxspace-abstract-smoother}{1} give, for the second term,
  \[
    \norm{\Dh}{\wh}^2 = (\Dh\wh, \wh)_{\dom} \overset{\eqref{eq:auxspace-abstract-smoother}_1}{\lesssim} c_{D,1} \norm{H^1(\dom)}{\wh}^2 + c_{D,0} \norm{L^2(\dom)}{\wh}^2 \overset{\eqref{eq:auxspace-abstract-decomposition}_{2,3}}{\lesssim} \left(c_{D,1} \frac{\H}{\h} + c_{D,0} \H^2\right) \norm{\ah}{\vh}^2.
  \]
  This gives Assumption \ref{ass:stable-decomposition} with $c_{0,H}^2 = 1$ and $c_{0,D}^2 = c_{D,1} \frac{\H}{\h} + c_{D,0} \H^2$. Applying Lemma \ref{lem:auxspace} concludes the proof.
\end{proof}

\section{Analysis tools}
\label{sec:analysis-tools}

We recall a series of well-known results that will be used throughout the paper. Let $K$ be a polytopal element of diameter $h_K$ that is star-shaped with respect to a ball of radius $\rho_K \geq c_R h_K$, as in Assumption~\ref{ass:coarse-mesh-regularity}(ii). In our setting, $K$ is either a cell of $\Th$ or an agglomerated element of $\TH$. These results are classical. We refer to \cite{Brenner-Scott-2008} for the standard simplicial setting, and to \cite{Cangiani2017hpVersion,DiPietro-Droniou-2020-HHO-Book} for their polytopal counterparts, where the hidden constants depend only on the space dimension $d$, the polynomial degree, and the shape regularity of $K$, that is the parameter $c_R$ of Assumption~\ref{ass:coarse-mesh-regularity} when $K \in \TH$ and the shape regularity of $\Th$ assumed in \sect{sec:preliminaries} when $K \in \Th$. The trace, Poincar\'e and Poincar\'e--Friedrichs inequalities, and the estimates for the elliptic projector with $m \leq 1$, rely on star-shapedness alone. The remaining projector estimates go through an inverse inequality on $K$, which for $K \in \TH$ uses in addition the simplicial covering of Assumption~\ref{ass:coarse-mesh-regularity}(iii).

\begin{itemize}
\item \emph{Scaled trace inequality:} For all $v \in H^1(K)$, it holds that
\begin{equation}
  \label{eq:trace-inequality-continuous}
  \norm{\partial K}{v}^2 \lesssim h_K^{-1} \norm{L^2(K)}{v}^2 + h_K \seminorm{1,K}{v}^2.
\end{equation}
\item \emph{Poincar\'e inequality:} For all $v \in H^1(K)\cap L^2_0(K)$, where $L^2_0(K)$ denotes the functions of $L^2(K)$ with zero mean on $K$, it holds that
\begin{equation}
  \label{eq:poincare-inequality}
  \norm{L^2(K)}{v} \lesssim h_K \seminorm{1,K}{v}.
\end{equation}
\item \emph{Poincar\'e--Friedrichs inequality:} For all $v \in H^1(K)$, with no constraint on $v$, it holds that
\begin{equation}
  \label{eq:poincare-friedrichs}
  \norm{L^2(K)}{v}^2 \lesssim h_K^2 \seminorm{1,K}{v}^2 + h_K \norm{\partial K}{v}^2.
\end{equation}
\item \emph{Approximation properties of local projectors:} We denote by $\lproj{K}{k} : L^2(K) \to \Pk{k}(K)$ the $L^2$-orthogonal projector on $K$ onto $\Pk{k}(K)$ polynomials. For all $v \in H^s(K)$ with $0 \leq s \leq k+1$, we have
  \begin{align}
    \seminorm{H^m(K)}{v - \lproj{K}{k} v} &\lesssim h_K^{s-m} \seminorm{H^s(K)}{v}, && 0 \leq m \leq s, \label{eq:l2-approx}\\
    \norm{\partial K}{v - \lproj{K}{k} v} &\lesssim h_K^{s-1/2} \seminorm{H^s(K)}{v}, && s \geq 1. \label{eq:l2-approx-trace}
  \end{align}
Similarly, denote by $\eproj{K}{k} : H^1(K) \to \Pk{k}(K)$ the elliptic projector defined by
\begin{equation}
  \label{eq:elliptic-projector}
  \begin{aligned}
    (\nabla \eproj{K}{k} v, \nabla w)_K &= (\nabla v, \nabla w)_K \quad \forall w \in \Pk{k}(K), \\
    (\eproj{K}{k} v, 1)_K &= (v, 1)_K.
  \end{aligned}
\end{equation}
For all $v \in H^s(K)$ with $1 \leq s \leq k+1$, we have
\begin{align}
  \seminorm{H^m(K)}{v - \eproj{K}{k} v} &\lesssim h_K^{s-m} \seminorm{H^s(K)}{v}, && 0 \leq m \leq s, \label{eq:h1-approx}\\
  \norm{\partial K}{v - \eproj{K}{k} v} &\lesssim h_K^{s-1/2} \seminorm{H^s(K)}{v}. \label{eq:h1-approx-trace}
\end{align}
\end{itemize}

Given a collection of polytopal elements $\mathcal{K}_h$ (e.g.\ $\Th$ or $\TH$), we denote by $\lproj{\mathcal{K}_h}{k}$ and $\eproj{\mathcal{K}_h}{k}$ the broken projectors acting element-wise, that is $\lproj{\mathcal{K}_h}{k} v|_K = \lproj{K}{k} (v|_K)$ and $\eproj{\mathcal{K}_h}{k} v|_K = \eproj{K}{k} (v|_K)$ for all $K \in \mathcal{K}_h$. The previous approximation properties hold element-wise, and thus also for the broken projectors.

\section{Additive Schwarz smoothers}
\label{sec:smoother}

Additive Schwarz smoothers have been widely studied in the literature in the context of \ac{dd} and \ac{gmg} solvers for conforming discretisations of elliptic problems. Classical examples include the \emph{point-Jacobi} smoother and overlapping block smoothers such as the \emph{vertex-star patch} smoother, in which each subdomain is the union of the cells sharing a mesh vertex. These are the two smoothers used in our experiments. We provide here a brief definition of these smoothers, and recall some of their properties that are relevant for the analysis of the \ac{gmg} method.

\subsection{Abstract setting}

Additive Schwarz smoothers are built from a decomposition of $\Uh$ into a family of subspaces. The family is indexed by a set $\mathcal{D}$ of geometrical entities of the mesh, such as nodes, vertices or patches of cells, the choice of $\mathcal{D}$ being what distinguishes one smoother from another. Let $\{W_r\}_{r \in \mathcal{D}}$ be subspaces of $\Uh$ such that $\Uh = \sum_{r \in \mathcal{D}} W_r$, and assume that the functions of $W_r$ are supported in the closure of an open set $\omega_r \subset \dom$.

We denote by $R_r : \Uh \to W_r$ the $L^2$-orthogonal projection onto $W_r$, and by $R_r^* : W_r \to \Uh$ its adjoint, which is the natural injection of $W_r$ into $\Uh$. We define the local bilinear forms as the restrictions of $\ah$ to the subspaces $W_r$, that is
\begin{equation}
  a_r(u_r,v_r) \doteq \ah(R_r^* u_r, R_r^* v_r) \quad \forall u_r, v_r \in W_r,
\end{equation}
with associated operators $A_r : W_r \to W_r$. The additive Schwarz smoother associated with the decomposition $\{W_r\}$ is then defined through its inverse as
\begin{equation}
  \label{eq:smoother-def-inverse}
  \Dh^{-1} \doteq \sum_{r \in \mathcal{D}} R_r^* A_r^{-1} R_r.
\end{equation}
The operator $\Dh$ is symmetric positive definite, and its quadratic form admits the classical characterisation (see \cite[Lemma 2.5]{Toselli-Widlund-DDM-book})
\begin{equation}
  \label{eq:smoother-def}
  (\Dh \vh, \vh)_{\dom} = \min \left\{ \sum_{r \in \mathcal{D}} a_r(v_r, v_r) \ : \ v_r \in W_r, \ \vh = \sum_{r \in \mathcal{D}} R_r^* v_r \right\}.
\end{equation}
We will refer to any set $\{v_r\}$, $v_r \in W_r$, such that $\vh = \sum_{r \in \mathcal{D}} R_r^* v_r$ as an \emph{admissible decomposition} of $\vh$. Any such decomposition provides an upper bound for $(\Dh\vh,\vh)_{\dom}$, which is how the first inequality in \eqref{eq:auxspace-abstract-smoother} is established in each of the cases below. The second inequality holds whenever the decomposition has finite overlap, in the following sense.

\begin{assumption} \label{ass:smoother-overlap}
  (see \cite[Assumption 3.2]{Toselli-Widlund-DDM-book}) The sets $\{\omega_r\}$ can be coloured using at most $\NO$ colours, in such a way that subregions with the same colour are disjoint. Since sets sharing a common point have different colours, every $x \in \Omega$ belongs to at most $\NO$ of them.
\end{assumption}

\begin{lemma}[Smoother stability]
  \label{lem:smoother-stability}
  Let Assumption \ref{ass:smoother-overlap} hold. Then Assumption~\ref{ass:smoother-stability} holds with $c_D^2 = \NO$, that is for all $\vh \in \Uh$,
  \[
    \ah(\vh, \vh) \leq \NO (\Dh \vh, \vh)_{\dom}.
  \]
\end{lemma}
\begin{proof}
  Let $\vh \in \Uh$ and let $\vh = \sum_{r \in \mathcal{D}} R_r^* v_r$, $v_r \in W_r$, be an arbitrary admissible decomposition. Grouping the subspaces by colour according to Assumption \ref{ass:smoother-overlap}, we can write $\vh = \sum_{c=1}^{\NO} w_c$ with $w_c \doteq \sum_{r \in \mathcal{C}_c} R_r^* v_r$, where $\mathcal{C}_c$ collects the entities of colour $c$. By the triangle and Cauchy-Schwarz inequalities,
  \[
    \ah(\vh, \vh) = \norm{\ah}{\textstyle\sum_{c} w_c}^2 \leq \Big(\sum_c \norm{\ah}{w_c}\Big)^2 \leq \NO \sum_{c} \norm{\ah}{w_c}^2,
  \]
  where the last step is the Cauchy--Schwarz inequality in $\mathbb{R}^{\NO}$. Since the sets $\omega_r$ of the same colour are disjoint, the supports of $\{R_r^* v_r\}_{r \in \mathcal{C}_c}$ meet only on sets of zero measure, and thus
  \[
    \norm{\ah}{w_c}^2 = \sum_{r \in \mathcal{C}_c} \ah(R_r^* v_r, R_r^* v_r) = \sum_{r \in \mathcal{C}_c} a_r(v_r, v_r).
  \]
  Combining both bounds gives $\ah(\vh,\vh) \leq \NO \sum_{r \in \mathcal{D}} a_r(v_r,v_r)$. Taking the minimum over all admissible decompositions and using the characterisation \eqref{eq:smoother-def} concludes the proof.
\end{proof}

For every decomposition considered below, $\NO$ is a fixed constant determined by the shape regularity of $\Th$ and by the polynomial degree. We therefore stop tracking it explicitly from here on, and absorb it into the hidden constants.

\subsection{Nodal decompositions}

The two smoothers used in our experiments correspond to two natural choices of the index set $\mathcal{D}$, and are covered by an elementary argument. Let $\Nhi$ denote the set of \emph{interior} Lagrangian nodes of $\Th$, that is those not lying on $\partial\dom$, and for $\nu \in \Nhi$ let $\varphi_\nu$ be the associated basis function. Since $\Uh \subset H^1_0(\dom)$, the family $\{\varphi_\nu\}_{\nu \in \Nhi}$ is a basis of $\Uh$ and
\begin{equation}
  \label{eq:nodal-decomposition}
  \vh = \sum_{\nu \in \Nhi} \vh(\nu) \varphi_\nu.
\end{equation}
We will use the standard scaling $\seminorm{1,\dom}{\varphi_\nu}^2 \simeq \h^{d-2}$ of the basis functions, together with the norm equivalence $\h^{d} \sum_{\nu \in \Nhi} |\vh(\nu)|^2 \simeq \norm{L^2(\dom)}{\vh}^2$ on $\Uh$.

\begin{corollary}[Point-Jacobi smoother]
  \label{cor:point-jacobi}
  Let $\mathcal{D} = \Nhi$ and let $\Dh$ be given by \eqref{eq:smoother-def-inverse} with $W_\nu \doteq \operatorname{span}\{\varphi_\nu\}$. Then \eqref{eq:auxspace-abstract-smoother} holds with $c_{D,0} \simeq \h^{-2}$, $c_{D,1} = 0$ and $c_D^2 \simeq 1$.
\end{corollary}
\begin{proof}
  Since $\{\varphi_\nu\}_{\nu \in \Nhi}$ is a basis of $\Uh$, the decomposition \eqref{eq:nodal-decomposition} is the only admissible one, and the minimum in \eqref{eq:smoother-def} is therefore attained at it. Hence
  \[
    (\Dh \vh, \vh)_{\dom} = \sum_{\nu \in \Nhi} |\vh(\nu)|^2 \seminorm{1,\dom}{\varphi_\nu}^2
    \simeq \h^{d-2} \sum_{\nu \in \Nhi} |\vh(\nu)|^2
    \simeq \h^{-2} \norm{L^2(\dom)}{\vh}^2,
  \]
  which is \eqp{eq:auxspace-abstract-smoother}{1} with $c_{D,0} \simeq \h^{-2}$ and $c_{D,1} = 0$. \Eqp{eq:auxspace-abstract-smoother}{2} follows from Lemma \ref{lem:smoother-stability}, since the supports of the $\varphi_\nu$ can be coloured with a number of colours bounded in terms of the shape regularity of $\Th$ and of $\k$.
\end{proof}

\begin{corollary}[Vertex-star patch smoother]
  \label{cor:vertex-patch}
  Let $\mathcal{D} = \mathcal{V}_h$ be the set of vertices of $\Th$ and, for $z \in \mathcal{V}_h$, let $S_z$ be its star, that is the union of the cells sharing $z$. Let $\Dh$ be given by \eqref{eq:smoother-def-inverse} with $W_z \doteq \{ \vh \in \Uh : \operatorname{supp}(\vh) \subset \ol{S_z} \}$. Then \eqref{eq:auxspace-abstract-smoother} holds with $c_{D,0} \simeq \h^{-2}$, $c_{D,1} = 0$ and $c_D^2 \simeq 1$.
\end{corollary}
\begin{proof}
  For each $\nu \in \Nhi$, let $e(\nu)$ be the lowest-dimensional mesh entity containing $\nu$, and fix one vertex $z(\nu)$ of $e(\nu)$. The basis function $\varphi_\nu$ is supported in the cells whose closure contains $\nu$, which are exactly those containing $e(\nu)$, and each of them contains the vertex $z(\nu)$ of $e(\nu)$. Hence $\operatorname{supp}(\varphi_\nu) \subset \ol{S_{z(\nu)}}$. Taking $e(\nu)$ of lowest dimension is essential here, since a vertex of a larger entity containing $\nu$ need not belong to every cell containing $\nu$. The map $\nu \mapsto z(\nu)$ partitions $\Nhi$, so
  \[
    v_z \doteq \sum_{\nu \,:\, z(\nu) = z} \vh(\nu)\varphi_\nu \in W_z
  \]
  regroups \eqref{eq:nodal-decomposition} into an admissible decomposition, since $\sum_{z \in \mathcal{V}_h} R_z^* v_z = \sum_{\nu \in \Nhi} \vh(\nu)\varphi_\nu = \vh$. The subspaces $W_z$ do overlap, but that plays no role here, because each node contributes to exactly one $v_z$.

  The number of nodes assigned to a given vertex is bounded in terms of the shape regularity of $\Th$ and of $\k$, so the triangle and Cauchy--Schwarz inequalities give $\seminorm{1,S_z}{v_z}^2 \lesssim \sum_{\nu \,:\, z(\nu) = z} |\vh(\nu)|^2 \seminorm{1,\dom}{\varphi_\nu}^2$. Summing over $z \in \mathcal{V}_h$, and using \eqref{eq:smoother-def} together with the fact that the map $\nu \mapsto z(\nu)$ is a partition,
  \[
    (\Dh \vh, \vh)_{\dom} \leq \sum_{z \in \mathcal{V}_h} \seminorm{1,S_z}{v_z}^2
    \lesssim \sum_{\nu \in \Nhi} |\vh(\nu)|^2 \seminorm{1,\dom}{\varphi_\nu}^2
    \simeq \h^{-2} \norm{L^2(\dom)}{\vh}^2.
  \]
  \Eqp{eq:auxspace-abstract-smoother}{2} follows again from Lemma \ref{lem:smoother-stability}. Here the overlap of the stars does matter, and it is controlled by Assumption \ref{ass:smoother-overlap}. Two stars intersect only if their vertices belong to a common cell, so the vertices can be coloured with a number of colours bounded by the shape regularity of $\Th$ in such a way that stars of the same colour are disjoint.
\end{proof}

\section{DG auxiliary space preconditioner}
\label{sec:dg}

We now instantiate the abstract framework of \sect{sec:auxspace}. The coarse auxiliary space $\UH$ is a broken polynomial space on the agglomerated mesh $\TH$, the coarse bilinear form $\aH$ is a Symmetric Interior Penalty Galerkin (SIPG) discretisation and the prolongation $\IHh = \IHhdg$ is built from an Oswald-type averaging operator on the fine mesh. We prove that the resulting preconditioned system has condition number bounded by $O(\H^3/\h^3)$  at fixed penalty, and by $O(\H^2/\h^2)$ for $\sigma \simeq \H/\h$, for the point-Jacobi and vertex-star patch smoothers of \sect{sec:smoother}.

\subsection{Coarse DG space and bilinear form}

Given a coarse polynomial degree $\K \geq 1$, in general independent of the fine degree $\k$, we define the coarse auxiliary space as the broken polynomial space of degree $\K$ on the agglomerated mesh:
\begin{equation}
  \label{eq:dg-space}
  \UH \doteq \Pk{\K}(\TH) = \bigl\{ \vH \in L^2(\dom) : \vH|_T \in \Pk{\K}(T) \;\forall T \in \TH \bigr\}.
\end{equation}
For an interior face $F \in \FHi$ shared by two coarse elements $T_1, T_2 \in \TH$, let $\bm{n}_{T_i}$ denote the unit normal to $F$ pointing outwards from $T_i$, and define the jump and the average as
\begin{equation}
  \label{eq:jump-average}
  \jump{\vH}|_F \doteq \vH|_{T_1} \bm{n}_{T_1} + \vH|_{T_2} \bm{n}_{T_2}, \qquad
  \avg{\nabla \vH}|_F \doteq \tfrac{1}{2}\bigl(\nabla \vH|_{T_1} + \nabla \vH|_{T_2}\bigr),
\end{equation}
where the restrictions denote traces taken from within each element.  For a boundary face $F \in \FHb$ with outward unit normal $\bm{n}_F$ we set $\jump{\vH}|_F \doteq \vH|_F \bm{n}_F$ and $\avg{\nabla \vH}|_F \doteq \nabla \vH|_F$, so that both notations extend uniformly to all $F \in \FH$. In either case $|\jump{\vH}|$ is the modulus of the difference of the two traces, respectively of the trace itself, so that $\norm{F}{\jump{\vH}}$ is the quantity appearing in all the estimates below. The same conventions are used on the fine mesh, on every $f \in \Fh$. The face size $H_F$ is as defined in \sect{sec:auxspace}, namely the mean diameter of the coarse elements sharing $F$ rather than the diameter of $F$ itself.

The coarse bilinear form $\aH : \UH \times \UH \to \mathbb{R}$ is taken as the SIPG formulation: for $\uH, \vH \in \UH$,
\begin{equation}
  \label{eq:sipg}
  \begin{aligned}
  \aH(\uH, \vH) &\doteq \sum_{T \in \TH} (\nabla \uH, \nabla \vH)_T 
     - \sum_{F \in \FH} \bigl[
        \bigl(\avg{\nabla \uH},\, \jump{\vH}\bigr)_F
        + \bigl(\avg{\nabla \vH},\, \jump{\uH}\bigr)_F
      \bigr] \\
    &\quad + \sum_{F \in \FH} \frac{\sigma}{H_F}\, \bigl(\jump{\uH},\, \jump{\vH}\bigr)_F,
  \end{aligned}
\end{equation}
where $\sigma > 0$ is the interior penalty parameter. The boundary face terms weakly enforce the homogeneous Dirichlet condition. The associated energy norm is
\begin{equation}
  \label{eq:dg-norm}
  \normDG{\vH}^2 \doteq \seminorm{1,\TH}{\vH}^2 + \sum_{F \in \FH} \frac{\sigma}{H_F} \norm{F}{\jump{\vH}}^2.
\end{equation}

The continuity and coercivity of the bilinear form \eqref{eq:sipg} with respect to $\normDG{\cdot}$ are well established in the literature, see, e.g.\ \cite{diPietro2012,Cangiani2017hpVersion}, and are given by the following lemma.
\begin{lemma} \label{lem:dg-norm-equivalence}
  If $\TH$ satisfies Assumption~\ref{ass:coarse-mesh-regularity}, then there is $\sigma_0 \geq 1$, depending only on $\K$ and the constants of Assumption~\ref{ass:coarse-mesh-regularity}, such that for all $\sigma \geq \sigma_0$
  \begin{align*}
    \aH(\uH, \vH) &\lesssim \normDG{\uH} \normDG{\vH} \ , & \forall \uH, \vH \in \UH. \\
    \aH(\vH, \vH) &\gtrsim \normDG{\vH}^2 \ , & \forall \vH \in \UH,
  \end{align*}
  with hidden constants independent of $\sigma$.
\end{lemma}

The $\sigma$-uniformity of the constants follows from Young's inequality together with the trace and inverse inequalities on $\TH$. We assume $\sigma \geq \sigma_0$ throughout, and track the dependence on $\sigma$ of all estimates below.
\subsection{Prolongation operator}

We define the prolongation $\IHhdg : \Pk{k}(\TH) \to \Uh$ by composing an element-wise interpolant with an Oswald-type nodal averaging operator, both defined on the fine mesh. Here $k \geq 1$ is a generic degree for the argument of $\IHhdg$, distinct from the fine degree $\k$ and the coarse degree $\K$. In this section $\IHhdg$ acts on the coarse \ac{dg} space and $k = \K$, while in \sect{sec:hho} it acts on the potential reconstruction of the coarse \ac{hho} function and $k = \K+1$. Since $\K$ is independent of $\k$, the argument of $\IHhdg$ may have a higher degree than the fine space can represent. With the choice $\K = \k$ used in our experiments, $k = \k$ here and $k = \k+1 > \k$ in \sect{sec:hho}. The estimates below assume no relation between $k$ and $\k$.

Let $\Ih$ denote the element-wise Lagrange interpolant of degree $\k$ onto $\Pk{\k}(\Th)$, i.e.\ $(\Ih w)|_t$ interpolates $w|_t$ at the Lagrangian nodes of $t$ for all $t \in \Th$. It is defined for any $w$ that is piecewise continuous in $\Th$, in particular on $\Pk{k}(\TH)$ for any $k \geq 1$, and it is the identity on $\Pk{\k}(\Th)$. For each interior node $\nu \in \Nhi$, let $\mathcal{N}_h(\nu) \subset \Th$ be the set of fine cells whose closure contains $\nu$, and let $\{\alpha_t(\nu)\}_{t \in \mathcal{N}_h(\nu)}$ be weights with $\alpha_t(\nu) \geq 0$ and $\sum_{t \in \mathcal{N}_h(\nu)} \alpha_t(\nu) = 1$. The \emph{averaging operator} $\Wh : \Pk{\k}(\Th) \to \Uh$ combines the (possibly multi-valued) values of its argument at each interior node with these weights and vanishes at boundary nodes:
\begin{equation}
  \label{eq:weighting-def}
  (\Wh \wh)(\nu) \doteq \sum_{t \in \mathcal{N}_h(\nu)} \alpha_t(\nu)\, (\wh|_t)(\nu) \ \ \text{ for } \nu \notin \partial\dom,
  \qquad (\Wh \wh)(\nu) \doteq 0 \ \ \text{ for } \nu \in \partial\dom.
\end{equation}
With uniform weights $\alpha_t(\nu) = |\mathcal{N}_h(\nu)|^{-1}$, $\Wh$ is the classical Oswald interpolation operator \cite{Karakashian2003,Burman2007}. The prolongation is the composition
\begin{equation}
  \label{eq:prolongation-def}
  \IHhdg \doteq \Wh \Ih : \Pk{k}(\TH) \to \Uh.
\end{equation}
Since $(\Ih \vH|_t)(\nu) = (\vH|_T)(\nu)$ for every fine cell $t \subset T$ and every node $\nu$ of $t$, the sum in \eqref{eq:weighting-def} can be grouped by agglomerates, and $\IHhdg$ reduces to a nodal averaging of the coarse one-sided values,
\begin{equation}
  \label{eq:prolongation-nodal}
  (\IHhdg\vH)(\nu) = \sum_{T \in \mathcal{N}_H(\nu)} \beta_T(\nu)\, (\vH|_T)(\nu) \ \ \text{ for } \nu \notin \partial\dom,
  \qquad \beta_T(\nu) \doteq \sum_{t \in \mathcal{N}_h(\nu),\, t \subset T} \alpha_t(\nu),
\end{equation}
where $\mathcal{N}_H(\nu) \subset \TH$ is the set of coarse elements whose closure contains $\nu$. The induced weights $\beta_T(\nu)$ are non-negative and sum to one. The resulting function is piecewise polynomial on $\Th$, globally continuous and vanishes on $\partial\dom$, so $\IHhdg\vH \in \Uh$.

\begin{remark}[Choice of the weights] \label{rem:weights}
  All the estimates below hold for any admissible choice of the weights $\alpha_t(\nu)$, with hidden constants independent of this choice. Two natural options are the uniform weights $\alpha_t(\nu) = |\mathcal{N}_h(\nu)|^{-1}$, for which $\beta_T(\nu)$ is proportional to the number of fine cells of $T$ at $\nu$, and the agglomerate-aware weights $\alpha_t(\nu) = \bigl(|\mathcal{N}_H(\nu)|\, n_T(\nu)\bigr)^{-1}$ for $t \subset T$, with $n_T(\nu)$ the number of fine cells of $T$ whose closure contains $\nu$, for which \eqref{eq:prolongation-nodal} is the plain average $\beta_T(\nu) = |\mathcal{N}_H(\nu)|^{-1}$ of the coarse one-sided values.
\end{remark}

The estimate we need for $\IHhdg$ rests on the following bound for the averaging operator $\Wh$. Classical versions of it, for the fine-mesh Oswald operator, are given in \cite[Th.~2.2]{Karakashian2003} for $r = 1$ and \cite[Lem.~3.2]{Burman2007} for $r = 0$. Neither is stated for general weights, and the $L^2$ bound is not stated for simplicial Lagrange elements with the boundary treatment used here, so we give the (elementary) proof.

\begin{lemma}[Nodal averaging on the fine mesh]
  \label{lem:nodal-averaging}
  For $r = 0, 1$ and all $\wh \in \Pk{\k}(\Th)$, the averaging operator $\Wh$ defined in \eqref{eq:weighting-def} with any weights that are non-negative and sum to one satisfies
  \begin{equation}
    \label{eq:nodal-averaging}
    \norm{H^r(\Th)}{\wh - \Wh \wh}^2 \lesssim \sum_{f \in \Fh} h_f^{1-2r} \norm{f}{\jump{\wh}}^2,
  \end{equation}
  with a hidden constant depending only on $\k$ and the shape regularity and quasi-uniformity of $\Th$, and in particular independent of the weights $\alpha_t$.
\end{lemma}
\begin{proof}
  Set $\delta_h \doteq \wh - \Wh \wh \in \Pk{\k}(\Th)$ and fix $t \in \Th$ and a node $\nu$ of $t$.

  If $\nu \notin \partial\dom$, then, since the weights sum to one,
  \[
    (\delta_h|_t)(\nu) = \sum_{t' \in \mathcal{N}_h(\nu)} \alpha_{t'}(\nu) \bigl[ (\wh|_t)(\nu) - (\wh|_{t'})(\nu) \bigr].
  \]
  Because $\Th$ is conforming and simplicial and $\dom$ is Lipschitz, the cells of $\mathcal{N}_h(\nu)$ tile a neighbourhood of $\nu$ in $\ol\dom$, so any two of them are connected by a chain
  \[
    t = t_0, t_1, \dots, t_M = t', \qquad t_i \in \mathcal{N}_h(\nu), \qquad M \leq |\mathcal{N}_h(\nu)|,
  \]
  in which $t_{i-1}$ and $t_i$ share a face $f_i \in \Fh(\nu)$, where $\Fh(\nu)$ denotes the set of fine faces containing $\nu$. Writing $\bm{n}_{i}$ for the unit normal to $f_i$ pointing outwards from $t_{i-1}$, the definition of the jump gives $(\wh|_{t_{i-1}})(\nu) - (\wh|_{t_i})(\nu) = \jump{\wh}|_{f_i}(\nu) \cdot \bm{n}_{i}$, so that each bracket telescopes as
  \[
    (\wh|_t)(\nu) - (\wh|_{t'})(\nu) = \sum_{i=1}^{M} \jump{\wh}|_{f_i}(\nu) \cdot \bm{n}_{i}.
  \]
  Since $|\mathcal{N}_h(\nu)|$ is bounded by shape regularity and $\alpha_{t'}(\nu) \leq 1$,
  \begin{equation}
    \label{eq:nodal-jump-bound}
    |(\delta_h|_t)(\nu)| \lesssim \max_{f \in \Fh(\nu)} |\jump{\wh}|_f(\nu)|.
  \end{equation}
If $\nu \in \partial\dom$, then $(\delta_h|_t)(\nu) = (\wh|_t)(\nu)$. Since $\partial\dom$ is a union of fine boundary faces, there is $f^{*} \in \Fhb$ with $\nu \in \ol{f^{*}}$. Let $t^{*} \in \mathcal{N}_h(\nu)$ be the cell of $f^{*}$. Chaining $t$ to $t^{*}$ as above, and using $\jump{\wh}|_{f^{*}} = \wh|_{f^{*}} \bm{n}_{f^{*}}$ on $\Fhb$,
  \[
    (\wh|_t)(\nu) = (\wh|_{t^{*}})(\nu) + \sum_{i} \jump{\wh}|_{f_i}(\nu) \cdot \bm{n}_{i}
    = \jump{\wh}|_{f^{*}}(\nu) \cdot \bm{n}_{f^{*}} + \sum_{i} \jump{\wh}|_{f_i}(\nu) \cdot \bm{n}_{i},
  \]
  so \eqref{eq:nodal-jump-bound} holds in this case as well. If $\nu$ is interior to $t$, then $\mathcal{N}_h(\nu) = \{t\}$, $\alpha_t(\nu) = 1$ and $(\delta_h|_t)(\nu) = 0$, so such nodes contribute nothing.

  Norm equivalence on the finite-dimensional space $\Pk{\k}(t)$, together with scaling, gives $\norm{L^2(t)}{\delta_h}^2 \simeq h_t^{d} \sum_{\nu} |(\delta_h|_t)(\nu)|^2$, the sum running over the nodes of $t$. The same equivalence on a face, applied to $\jump{\wh}|_f \in \Pk{\k}(f)$, gives $\norm{f}{\jump{\wh}}^2 \simeq h_f^{d-1} \sum_{\nu} |\jump{\wh}|_f(\nu)|^2$. Combining these with \eqref{eq:nodal-jump-bound}, summing over $t \in \Th$, and using that each face is counted a bounded number of times together with the quasi-uniformity $h_t \simeq h_f \simeq \h$,
  \[
    \norm{L^2(\Th)}{\delta_h}^2 \lesssim \h^{d} \sum_{f \in \Fh} \sum_{\nu \in f} |\jump{\wh}|_f(\nu)|^2
    \simeq \h^{d} \sum_{f \in \Fh} h_f^{-(d-1)} \norm{f}{\jump{\wh}}^2
    \simeq \sum_{f \in \Fh} h_f \norm{f}{\jump{\wh}}^2,
  \]
  which is \eqref{eq:nodal-averaging} for $r = 0$. The case $r = 1$ follows from the case $r = 0$ and the inverse inequality $\seminorm{1,t}{\delta_h} \lesssim h_t^{-1} \norm{L^2(t)}{\delta_h}$ on each fine cell, using $h_t \simeq h_f \simeq \h$ and $h_f \lesssim \mathrm{diam}(\dom)$ for the $L^2$ part.
\end{proof}
\begin{lemma}[Bounds for $\Wh\Ih$ on the coarse broken space]
  \label{lem:oswald}
  Let $k \geq 1$ and $\vH \in \Pk{k}(\TH)$. The operator $\IHhdg = \Wh\Ih$ satisfies, for $r = 0, 1$,
  \begin{equation}
    \norm{H^r(\TH)}{\vH - \IHhdg\vH}^2 \lesssim \sum_{F \in \FH} \h^{1-2r} \norm{F}{\jump{\vH}}^2 + \sum_{T \in \TH} \left(\frac{\h}{H_T}\right)^{2(\k+1-r)} H_T^{2(1-r)} \seminorm{1,T}{\vH}^2, \label{eq:oswald}
  \end{equation}
  The hidden constant depends only on $k$, $\k$, the shape regularity and quasi-uniformity of $\Th$ and Assumption~\ref{ass:coarse-mesh-regularity}, and is independent of the weights of $\Wh$.
\end{lemma}
\begin{proof}
  Set $\wh \doteq \Ih \vH \in \Pk{\k}(\Th)$, so that $\IHhdg \vH = \Wh \wh$.
  Since $\vH - \IHhdg\vH$ belongs to $H^r(T)$ on every $T \in \TH$ (both $\vH|_T$ and the globally continuous $\IHhdg\vH$ do), its $\TH$- and $\Th$-broken norms coincide, and we may work with the latter. Since $\wh$ is a priori only piecewise smooth on $\Th$, the intermediate terms below are $\Th$-broken norms. We can thus split
  \[
    \norm{H^r(\TH)}{\vH - \IHhdg\vH}^2 = \norm{H^r(\Th)}{\vH - \IHhdg\vH}^2 \lesssim \norm{H^r(\Th)}{\vH - \wh}^2 + \norm{H^r(\Th)}{\wh - \Wh \wh}^2 \doteq \mathfrak{I}_1 + \mathfrak{I}_2.
  \]
  
  For $\mathfrak{I}_1$, standard interpolation estimates on each fine cell $t \subset T$ give
  \[
    \mathfrak{I}_1 = \sum_{T \in \TH} \sum_{t \subset T} \norm{H^r(t)}{\vH - \Ih \vH}^2
    \lesssim \sum_{T \in \TH} \h^{2(\k+1-r)} \seminorm{H^{\k+1}(T)}{\vH}^2
    \lesssim \sum_{T \in \TH} \left(\frac{\h}{H_T}\right)^{2(\k+1-r)} H_T^{2(1-r)} \seminorm{1,T}{\vH}^2,
  \]
  where the last step uses the inverse inequality $\seminorm{H^{\k+1}(T)}{\vH} \lesssim H_T^{-\k} \seminorm{1,T}{\vH}$, valid for polynomials of fixed degree $k$ on elements satisfying Assumption~\ref{ass:coarse-mesh-regularity} \cite{Cangiani2017hpVersion}. If $k \leq \k$, then $\Ih \vH = \vH$ and $\mathfrak{I}_1 = 0$.

  We now bound $\mathfrak{I}_2$. Lemma \ref{lem:nodal-averaging} gives
  \[
    \mathfrak{I}_2 \lesssim \sum_{f \in \Fh} h_f^{1-2r} \norm{f}{\jump{\wh}}^2.
  \]
  For simplicial Lagrange elements, the trace of $(\Ih w)|_t$ on a face $f \in \Fh(t)$ is determined by the values of $w|_t$ at the Lagrangian nodes lying on $f$. As a consequence, for a fine face $f$ shared by two cells $t_1, t_2$ contained in the same agglomerate $T$, the traces of $\wh$ from both sides interpolate the same polynomial $\vH|_T$ at the same nodes, and thus $\jump{\wh}|_f = 0$. For a fine face $f \subset F \in \FH$, the same argument shows that $\jump{\wh}|_f = \Ih^f (\jump{\vH}|_f)$, where $\Ih^f$ denotes the Lagrange interpolant of degree $\k$ on $f$. By the $L^2(f)$-stability of $\Ih^f$ on the finite-dimensional space $\Pk{k}(f)$, we conclude that $\norm{f}{\jump{\wh}} \lesssim \norm{f}{\jump{\vH}}$, and thus
  \[
    \mathfrak{I}_2 \lesssim \sum_{F \in \FH} \sum_{f \subset F} h_f^{1-2r} \norm{f}{\jump{\vH}}^2 \lesssim \sum_{F \in \FH} \h^{1-2r} \norm{F}{\jump{\vH}}^2,
  \]
  where the last step uses that the fine faces $f \subset F$ tile $F$ and that $h_f \simeq \h$ for every $f \in \Fh$, by the quasi-uniformity of $\Th$ assumed in \sect{sec:preliminaries}, with a constant independent of $\h$ and $\H$.
  Combining the bounds for $\mathfrak{I}_1$ and $\mathfrak{I}_2$ concludes the proof.
\end{proof}

\begin{lemma} \label{lem:oswald-dg}
  For any $k \geq 1$ and $\vH \in \Pk{k}(\TH)$, the operator $\IHhdg = \Wh\Ih$ satisfies
  \begin{equation}
    \norm{H^1(\dom)}{\IHhdg \vH}^2 \lesssim \Bigl(1 + \frac{\H}{\sigma\h}\Bigr) \normDG{\vH}^2.
  \end{equation}
  In particular, since $\sigma \geq \sigma_0 \geq 1$ and $\h \leq \H$, the coarser bound $\norm{H^1(\dom)}{\IHhdg \vH}^2 \lesssim \frac{\H}{\h} \normDG{\vH}^2$ also holds.
\end{lemma}
\begin{proof}
  Since $\IHhdg\vH \in H^1_0(\dom)$, a Poincar\'e inequality gives $\norm{H^1(\dom)}{\IHhdg \vH} \lesssim \seminorm{1,\dom}{\IHhdg \vH}$. By the triangle inequality,
  \[
    \seminorm{1,\dom}{\IHhdg \vH} \leq \seminorm{1,\TH}{\vH - \IHhdg \vH} + \seminorm{1,\TH}{\vH}.
  \]
  By Lemma \ref{lem:oswald} with $r = 1$ and $\h \leq H_T$,
  \[
    \seminorm{1,\TH}{\vH - \IHhdg \vH}^2 \lesssim \sum_{F \in \FH} \h^{-1} \norm{F}{\jump{\vH}}^2 + \seminorm{1,\TH}{\vH}^2.
  \]
  Using the fact that $H_F \simeq \H$, the jump term can be bounded as
  \[
    \sum_{F \in \FH} \h^{-1} \norm{F}{\jump{\vH}}^2 \lesssim \sigma^{-1}\frac{\H}{\h} \sum_{F \in \FH} \frac{\sigma}{H_F} \norm{F}{\jump{\vH}}^2 \leq \sigma^{-1}\frac{\H}{\h} \normDG{\vH}^2.
  \]
  Since $\seminorm{1,\TH}{\vH}^2 \leq \normDG{\vH}^2$, combining the two bounds concludes the proof.
\end{proof}

\begin{lemma} \label{lem:oswald-prolongation}
  Assume that, given $\vh \in \Uh$, there exists $\vH \in \Pk{k}(\TH)$, $k \geq 1$, such that for $r = 0, 1$ it holds that
  \begin{equation}
    \label{eq:coarse-approx}
    \norm{H^r(\TH)}{\vh - \vH} \lesssim \H^{1-r} \seminorm{1,\dom}{\vh}.
  \end{equation}
  Then the operator $\IHhdg = \Wh\Ih$ satisfies
  \begin{equation}
    \label{eq:prolongation-approx}
    \norm{L^2(\dom)}{\vh - \IHhdg\vH}^2 \lesssim \H^2 \seminorm{1,\dom}{\vh}^2 \quad \text{and} \quad
    \norm{H^1(\dom)}{\vh - \IHhdg\vH}^2 \lesssim \frac{\H}{\h} \seminorm{1,\dom}{\vh}^2,
  \end{equation}
  and, if $\vH \in \UH$,
  \begin{equation}
    \label{eq:dg-bound}
    \normDG{\vH} \lesssim (1+\sigma)^{1/2}\, \seminorm{1,\dom}{\vh}.
  \end{equation}
\end{lemma}
\begin{proof}
  We first show the auxiliary bounds
  \begin{equation}
    \label{eq:coarse-approx-core}
    \seminorm{1,\TH}{\vH}^2 \lesssim \seminorm{1,\dom}{\vh}^2
    \quad \text{and} \quad
    \sum_{F \in \FH} H_F^{-1} \norm{F}{\jump{\vH}}^2 \lesssim \seminorm{1,\dom}{\vh}^2.
  \end{equation}
  Using the triangle inequality, together with \eqref{eq:coarse-approx}, we can bound
  \[
    \seminorm{1,\TH}{\vH}^2 \lesssim \seminorm{1,\TH}{\vh}^2 + \seminorm{1,\TH}{\vh - \vH}^2 \lesssim \seminorm{1,\dom}{\vh}^2.
  \]
  For the jump term, we can use the fact that $\vh$ is continuous across interior $\FH$ faces and vanishes on the boundary ones to write
  \[
    \sum_{F \in \FH} H_F^{-1} \norm{F}{\jump{\vH}}^2 = \sum_{F \in \FH} H_F^{-1} \norm{F}{\jump{\vh - \vH}}^2.
  \]
  Then noting that $H_F \simeq \H$, rearranging the face terms into element terms and using the scaled trace inequality \eqref{eq:trace-inequality-continuous} on each subdomain $T$, we have
  \[
    \sum_{F \in \FH} H_F^{-1} \norm{F}{\jump{\vH}}^2 \lesssim \sum_{T \in \TH} \H^{-1} \norm{\partial T}{\vh - \vH}^2 \lesssim \sum_{T \in \TH} \H^{-1} \left( H_T^{-1} \norm{L^2(T)}{\vh - \vH}^2 + H_T \norm{H^1(T)}{\vh - \vH}^2\right).
  \]
  Using the fact that $H_T \simeq \H$ for all $T \in \TH$, summing over $T \in \TH$ and using \eqref{eq:coarse-approx}, we obtain
  \[
    \sum_{F \in \FH} H_F^{-1} \norm{F}{\jump{\vH}}^2 \lesssim \H^{-1} \left( \H^{-1} \norm{L^2(\dom)}{\vh - \vH}^2 + \H \norm{H^1(\TH)}{\vh - \vH}^2 \right) \lesssim \seminorm{1,\dom}{\vh}^2,
  \]
  which proves \eqref{eq:coarse-approx-core}.

  We now prove \eqref{eq:prolongation-approx}. For $r = 0, 1$, a triangle inequality gives
  \[
    \norm{H^r(\dom)}{\vh - \IHhdg\vH}^2 \lesssim \norm{H^r(\TH)}{\vh - \vH}^2 + \norm{H^r(\TH)}{\vH - \IHhdg\vH}^2 \doteq \mathfrak{T}_1 + \mathfrak{T}_2.
  \]
  By \eqref{eq:coarse-approx}, $\mathfrak{T}_1 \lesssim \H^{2(1-r)} \seminorm{1,\dom}{\vh}^2$. For $\mathfrak{T}_2$, we can use Lemma \ref{lem:oswald} and $\h \leq H_T$, together with $H_F \simeq \H$ and \eqref{eq:coarse-approx-core}, to write
  \[
    \begin{aligned}
      \mathfrak{T}_2 &\lesssim \sum_{F \in \FH} \h^{1-2r} \norm{F}{\jump{\vH}}^2 + \H^{2(1-r)} \seminorm{1,\TH}{\vH}^2
      \lesssim \H \h^{1-2r} \sum_{F \in \FH} H_F^{-1} \norm{F}{\jump{\vH}}^2 + \H^{2(1-r)} \seminorm{1,\TH}{\vH}^2 \\
      &\lesssim \bigl(\H \h^{1-2r} + \H^{2(1-r)}\bigr) \seminorm{1,\dom}{\vh}^2.
    \end{aligned}
  \]
  Since $\h \leq \H$, the factor in brackets is bounded by $\H^2$ for $r = 0$ and by $\H/\h$ for $r = 1$, which gives \eqref{eq:prolongation-approx}.

  Finally, if $\vH \in \UH$, the definition \eqref{eq:dg-norm} of the \ac{dg} norm and \eqref{eq:coarse-approx-core} give
  \[
    \normDG{\vH}^2 = \seminorm{1,\TH}{\vH}^2 + \sigma \sum_{F \in \FH} H_F^{-1} \norm{F}{\jump{\vH}}^2 \lesssim (1 + \sigma) \seminorm{1,\dom}{\vh}^2,
  \]
  which is \eqref{eq:dg-bound}.
\end{proof}

\subsection{Main result}

\begin{theorem}
  \label{th:dg-main}
  The \ac{dg} discretisation, with the smoother described in Section \ref{sec:smoother} and the prolongation $\IHh = \IHhdg$ of \eqref{eq:prolongation-def}, fulfils the assumptions in Lemma \ref{lem:auxspace-concrete}. More precisely, the prolongation stability holds with $c_I^2 \lesssim \max(1, \H/(\sigma\h))$ and the first bound of the stable decomposition with $c_{0,H}^2 \lesssim 1 + \sigma$. All remaining hidden constants are independent of $\sigma$.
\end{theorem}
\begin{proof}
  That the smoother satisfies Equation \eqref{eq:auxspace-abstract-smoother} is given by Corollary \ref{cor:point-jacobi} or \ref{cor:vertex-patch}, depending on the smoother at hand. \\
  Let $\vH \in \UH$. Using Lemma \ref{lem:dg-norm-equivalence} together with Lemma \ref{lem:oswald-dg}, we obtain
  \[
    \norm{\ah}{\IHhdg\vH}^2 
    = \seminorm{1,\dom}{\IHhdg\vH}^2 
    \lesssim \Bigl(1 + \frac{\H}{\sigma\h}\Bigr) \normDG{\vH}^2 
    \lesssim \Bigl(1 + \frac{\H}{\sigma\h}\Bigr) \norm{a_H}{\vH}^2
    \lesssim \frac{\H}{\h} \norm{a_H}{\vH}^2
  \]
  which proves Equation \eqref{eq:auxspace-abstract-prolongation}.
  
  We now turn to Equation \eqref{eq:auxspace-abstract-decomposition}. Given $\vh \in \Uh$, we choose as coarse approximant $\vH \doteq \lproj{\TH}{\K} \vh \in \UH$. From standard approximation properties of $\lproj{\TH}{\K}$ \eqref{eq:l2-approx}, applied with $s = 1$ and $m = 0, 1$, we have that for $r = 0, 1$
  \[
    \norm{H^r(\TH)}{\vh - \vH} \lesssim \H^{1-r} \seminorm{1,\dom}{\vh},
  \]
  and thus our coarse approximant $\vH$ fulfils the assumption \eqref{eq:coarse-approx} of Lemma \ref{lem:oswald-prolongation}.

  Using \eqref{eq:prolongation-approx}, together with $\seminorm{1,\dom}{\vh} = \norm{\ah}{\vh}$, we obtain
  \[
    \norm{L^2(\dom)}{\vh - \IHhdg \vH}^2 \lesssim \H^2 \norm{\ah}{\vh}^2 \quad \text{ and } \quad \norm{H^1(\dom)}{\vh - \IHhdg \vH}^2 \lesssim \frac{\H}{\h} \norm{\ah}{\vh}^2
  \]
  Similarly, using \eqref{eq:dg-bound} together with Lemma \ref{lem:dg-norm-equivalence}, we can write
  \[
    \norm{a_H}{\vH} \lesssim \normDG{\vH} \lesssim (1+\sigma)^{1/2}\, \seminorm{1,\dom}{\vh} = (1+\sigma)^{1/2}\, \norm{\ah}{\vh}
  \]
  which proves \eqp{eq:auxspace-abstract-decomposition}{2,3} and \eqp{eq:auxspace-abstract-decomposition}{1}, respectively.
\end{proof}

\begin{remark} \label{rem:dg-final}
  Theorem \ref{th:dg-main} supplies $c_I^2 \lesssim \max(1, \H/(\sigma\h))$, $c_{0,H}^2 \lesssim 1 + \sigma$ and $c_{0,D}^2 \lesssim c_{D,1} \frac{\H}{\h} + c_{D,0} \H^2$, and inserting these into Lemma \ref{lem:auxspace} gives the final bound
  \[
    \kappa(P_{\mathrm{add,DG}}^{-1} \Ah) \lesssim \Bigl(c_D^2 + \lambda_1 \max\Bigl(1, \frac{\H}{\sigma\h}\Bigr)\Bigr)\Bigl(\frac{\H^2}{\h^2} + \frac{1+\sigma}{\lambda_0}\Bigr),
  \]
  where we have used that Corollaries \ref{cor:point-jacobi} and \ref{cor:vertex-patch} give $c_{D,0} \simeq \h^{-2}$ and $c_{D,1} = 0$ for both smoothers considered, and that $c_D^2 \simeq 1$. The coarse solvers of \sect{sec:multilevel} contribute the factors $\lambda_0$ and $\lambda_1$. For an exact coarse solve, $\lambda_0 = \lambda_1 = 1$, and two regimes of the penalty are of interest. For fixed $\sigma = O(1)$ the bound is $O(\H^3/\h^3)$, in agreement with similar bounds found for traditional \ac{fe} methods. Choosing instead the penalty proportional to the coarsening ratio, $\sigma \simeq \H/\h$, makes the first factor $O(1)$ while the second remains $O(\H^2/\h^2)$, so that
  \[
    \kappa(P_{\mathrm{add,DG}}^{-1} \Ah) \lesssim \frac{\H^2}{\h^2},
  \]
  one power of $\H/\h$ better.  \sect{sec:multilevel} shows how this bound composes when the coarse polytopal problem is itself solved approximately. See \cite{Chan1994} or the introduction to \cite{Dryja2016} for a discussion on the dependence of non-overlapping additive Schwarz subdomain sizes on the final bounds.
\end{remark}

\section{HHO auxiliary space preconditioner}
\label{sec:hho}

We now instantiate the abstract framework of \sect{sec:auxspace} with a coarse space built from a Hybrid High-Order (HHO) discretisation on $\TH$. Specifically, the coarse auxiliary space will be the skeleton space obtained after elimination of the cell unknowns via static condensation. The HHO bilinear form on the skeleton controls the broken $H^1$ seminorm and the interface jumps of the associated potential reconstruction, which reduces the analysis to the tools of \sect{sec:dg}.

\subsection{HHO discretisation on the coarse mesh}

We introduce an HHO discretisation on the coarse mesh $\TH$, using a mixed-order approach that employs polynomials of degree $\K+1$ for the cell unknowns and polynomials of degree $\K$ for the face unknowns. This choice admits a simpler stabilisation term than the original \ac{hho} method \cite{DiPietro-Droniou-2020-HHO-Book}, not involving the potential reconstruction, while still providing optimal approximation properties. After static condensation, the size of the global skeleton system is the same as for the equal-order method.

\emph{Face polynomial spaces.}
Given a face $F \in \FH$, we denote our skeleton space on $F$ by $\IPk{\K}(F)$. Recall that $F$ may be non-planar, but rather a collection of planar faces from $\Fh$. The original \ac{hho} method \cite{DiPietro-Droniou-2020-HHO-Book} uses the Cartesian product of polynomial spaces on each planar face $f \in \Fh$ to define the face unknowns, that is $\IPk{\K}(F) = \bigtimes_{f \in \Fh(F)} \Pk{\K}(f)$, where, extending the notation of \sect{sec:preliminaries}, $\Fh(F) \doteq \{ f \in \Fh : f \subset F\}$ denotes the set of fine faces composing $F$. However, this choice leads to a large number of degrees of freedom on non-planar faces. Instead, we consider the face polynomial space introduced in \cite{Yemm2022ANA}:
\begin{equation}
  \label{eq:face-space}
  \IPk{\K}(F) \doteq \Pk{0}(\dom)|_F + \nabla \Pk{\K+1}(\dom)|_F \cdot \bm{n}_F,
\end{equation}
where $\bm{n}_F$ is a (piecewise-constant) unit normal to $F$, either of the two on an interior face, the space being unchanged by the choice. For planar faces $\IPk{\K}(F)$ reduces to the standard choice $\Pk{\K}(F)$, but for non-planar agglomerated interfaces this choice provides optimal approximation properties with fewer degrees of freedom. These spaces have been shown to be effective in the context of agglomeration-based methods in \cite{badia2026geometricmultigridsolvershybrid}. We denote by $\IPk{\K}(\FH) = \{ (\lambda_F)_F : \lambda_F \in \IPk{\K}(F) , \ F \in \FH \}$ its corresponding broken skeleton space. We define $\lproj{F}{\K}$ as the $L^2$-orthogonal projector onto $\IPk{\K}(F)$, and $\lproj{\partial T}{\K}$ as the broken projector acting face-wise on $\partial T$, that is $\lproj{\partial T}{\K} v|_F = \lproj{F}{\K}(v|_F)$ for all $F \in \FH(T)$. Similarly, we define $\lproj{\FH}{\K}$ as the broken projector acting face-wise on $\FH$. Throughout this section, projectors carrying a face subscript ($F$, $\partial T$ or $\FH$) always denote the $L^2$-orthogonal projections onto the interface spaces $\IPk{\K}(\cdot)$ just defined, while projectors carrying a cell or fine-mesh subscript denote the projections onto full polynomial spaces of Section \ref{sec:analysis-tools}. Both conventions coexist in the formulas below.
\emph{Hybrid formulation.} The local HHO space of degree $\K$ on $T$ is
\[
  \UHhyb(T) \doteq \Pk{\K+1}(T) \times M_{\partial T} \ , \quad M_{\partial T} \doteq \prod_{F \in \FH(T)} \IPk{\K}(F),
\]
with hybrid unknowns $\vThyb = (v_T, v_{\partial T})$, where $v_T \in \Pk{\K+1}(T)$ is the cell unknown and $v_{\partial T} = (v_F)_{F \in \FH(T)}$ collects the face unknowns. The global HHO space is
\begin{equation}
  \label{eq:hho-space}
  \UHhyb \doteq U_{\TH} \times \MH \doteq \Pk{\K+1}(\TH) \times \IPk{\K}_0(\FH) \ , \quad \IPk{\K}_0(\FH) \doteq \bigl\{ \lH \in \IPk{\K}(\FH) \ : \ \lambda_F = 0 \ \ \forall F \in \FHb \bigr\},
\end{equation}
with global hybrid unknowns $\vHhyb = (v_{\TH}, v_{\FH})$ and $\vHhyb|_T = \vThyb$. The constraint on boundary faces strongly enforces the homogeneous Dirichlet boundary condition on the skeleton. Together with the norm equivalence of Lemma \ref{lem:hho-norm-equivalence}, it ensures that the condensed bilinear form $\aH$ introduced below is positive definite on $\MH$, see Lemma \ref{lem:hho-pd} (without the boundary constraint, the constant skeleton function would belong to its kernel).
On the hybrid space $\UHhyb$ we define the usual hybrid seminorms (the local one has the one-dimensional kernel $\{(c, c|_{\partial T})\}$, since constants belong to the face spaces, while the global one is a norm on $\UHhyb$ by the positive-definiteness argument below) 
\begin{equation}
  \label{eq:hho-hybrid-norm}
  \normHHO{\vHhyb}^2 \doteq \sum_{T \in \TH} \normHHOloc{\vThyb}^2 
  \ , \quad
  \normHHOloc{\vThyb}^2 \doteq \seminorm{1,T}{v_{T}}^2 + \HT^{-1} \norm{\partial T}{v_{T} - v_{\partial T}}^2.
\end{equation}

The local potential reconstruction $\rec{T}{\K+1} : \UHhyb(T) \to \Pk{\K+1}(T)$ is the unique polynomial satisfying, for all $w \in \Pk{\K+1}(T)$,
\begin{equation}
  \label{eq:hho-rec}
  \begin{aligned}
  (\nabla \rec{T}{\K+1}\vThyb, \nabla w)_T &= (\nabla v_T, \nabla w)_T + (v_{\partial T} - v_T,\, \bm{n}_{\partial T} \cdot \nabla w)_{\partial T}, \\
  (\rec{T}{\K+1}\vThyb, 1)_T &= (v_T, 1)_T,
  \end{aligned}
\end{equation}
where $\bm{n}_{\partial T}$ is the unit outward normal to $T$. The global reconstruction $\rec{H}{\K+1} : \UHhyb \to \Pk{\K+1}(\TH)$ applies $\rec{T}{\K+1}$ element-wise.

\noindent The local hybrid bilinear form on $T$ reads
\begin{equation}
  \label{eq:hho-form-local}
  \ul{a}_{T}(\uThyb, \vThyb) \doteq
  (\nabla \rec{T}{\K+1}\uThyb, \nabla \rec{T}{\K+1}\vThyb)_T
  + s_T(\uThyb, \vThyb),
\end{equation}
with stabilisation
\begin{equation}
  \label{eq:hho-stab}
  s_T(\uThyb, \vThyb) \doteq \HT^{-1} \left(\lproj{\partial T}{\K}u_T - u_{\partial T},\, \lproj{\partial T}{\K}v_T - v_{\partial T}\right)_{\partial T}.
\end{equation}
The global hybrid bilinear form is then $\ul{a}_{H}(\uHhyb, \vHhyb) \doteq \sum_{T \in \TH} \ul{a}_{T}(\uThyb, \vThyb)$. Unlike the SIPG method of Section \ref{sec:dg}, whose penalty parameter must be taken large enough for coercivity and is later exploited at $\sigma \simeq \H/\h$, the stabilisation \eqref{eq:hho-stab} carries no penalty coefficient. With the unit weight, the local bilinear form is already equivalent to the hybrid seminorm (Lemma \ref{lem:hho-norm-equivalence}) with constants of order one, and no super-penalisation is used in this section.

\emph{Condensed formulation.} The cell unknowns can be locally eliminated by exploiting the block structure of the hybrid bilinear form. For each $T \in \TH$, we define the inverse static condensation operator $\isc{T} : M_{\partial T} \to \Pk{\K+1}(T)$ as the unique solution of
\begin{equation}
  \label{eq:inverse-static-condensation}
  \ul{a}_{T}\bigl((\isc{T}\lH, 0),\, (w_T, 0)\bigr)
  = -\ul{a}_{T}\bigl((0, \lH),\, (w_T, 0)\bigr)
  \quad \forall w_T \in \Pk{\K+1}(T).
\end{equation}
That is, $\isc{T}\lH$ recovers the cell unknown from the skeleton data $\lH|_{\partial T}$ by solving the local cell equation. The local problem \eqref{eq:inverse-static-condensation} is well-posed, since $\ul{a}_{T}$ is positive definite on cell bubbles $(w_T, 0)$. For the face spaces \eqref{eq:face-space}, this follows from the inclusion $\nabla \Pk{\K+1}(T) \cdot \bm{n}_F|_F \subset \IPk{\K}(F)$ together with $\Pk{0}(\dom)|_F \subset \IPk{\K}(F)$, which is what rules out the constants and is why \eqref{eq:face-space} carries the $\Pk{0}$ term on non-planar faces, see \cite{Yemm2022ANA,badia2026geometricmultigridsolvershybrid}. The global broken operator $\iscH : \MH \to \Pk{\K+1}(\TH)$ acts element-wise: $(\iscH \lH)|_T \doteq \isc{T}(\lH|_{\partial T})$ for all $T \in \TH$. We also define the hybrid extension $\iscHhyb : \MH \to \UHhyb$ by
\begin{equation}
  \label{eq:hybrid-extension}
  \iscHhyb \lH \doteq (\iscH \lH,\, \lH).
\end{equation}

Substituting the inverse static condensation into the hybrid bilinear form yields the local and global skeleton bilinear forms:
\begin{align}
  a_{T}(\lH, \mH) &\doteq \ul{a}_{T}\bigl((\isc{T}\lH, \lH),\, (\isc{T}\mH, \mH)\bigr), \label{eq:skeleton-bilinear-form-local} \\
  \aH(\lH, \mH) &\doteq \sum_{T \in \TH} a_{T}(\lH, \mH). \label{eq:skeleton-bilinear-form-global}
\end{align}
We identify the coarse auxiliary space as $\UH \doteq \MH$ equipped with the skeleton bilinear form $\aH$, and denote by $\AH : \MH \to \MH$ the associated operator, characterised by $(\AH\lH, \mH) = \aH(\lH, \mH)$, where $(\cdot,\cdot)$ is the $L^2(\FH)$ inner product.

Since $\iscH\lH$ minimises the local hybrid energy over all cell unknowns with the same skeleton data $\lH$, the skeleton bilinear form satisfies the energy minimisation property: for any $\lH \in \MH$ and any $\vHhyb = (v_H, \lH) \in \UHhyb$ with skeleton values $\lH$,
\begin{equation}
  \label{eq:energy-minimisation}
  \aH(\lH, \lH) \leq \ul{a}_H(\vHhyb, \vHhyb).
\end{equation}

\subsection{HHO analysis tools}

We start by recalling the following bound, whose proof can be found in \cite[Lemma 7]{Yemm2022ANA}: for all $v \in H^1(T)$,
\begin{equation} \label{eq:hho-intermediate-trace-bound-1}
  \seminorm{1,T}{\lproj{T}{\K}v}^2 + \HT^{-1} \norm{\partial T}{\lproj{\partial T}{\K}v - \lproj{T}{\K} v}^2 \lesssim \seminorm{1,T}{v}^2.
\end{equation}
We note that the proof of \eqref{eq:hho-intermediate-trace-bound-1} relies on the inclusion $\Pk{0}(\dom)|_F \subset \IPk{\K}(F)$, which holds by construction of \eqref{eq:face-space}.
From \eqref{eq:hho-intermediate-trace-bound-1} we can further deduce the following result
\begin{equation} \label{eq:hho-intermediate-trace-bound-2}
  \HT^{-1} \norm{\partial T}{\lproj{\partial T}{\K}v - v}^2 \lesssim \seminorm{1,T}{v}^2, 
\end{equation}
where it suffices to add and subtract $\lproj{T}{\K}v$ in the face term and use the triangle inequality followed by \eqref{eq:hho-intermediate-trace-bound-1} and \eqref{eq:l2-approx-trace}.

We now prove the following two bounds.
\begin{lemma} \label{lem:hho-local-bounds}
  For all $T \in \TH$ and $\vThyb \in \UHhyb(T)$,
  \begin{equation} \label{eq:hho-cell-coercivity}
    \seminorm{1,T}{v_T}^2 + \HT^{-1} \norm{\partial T}{v_{\partial T} - \lproj{\partial T}{\K} v_T}^2 \lesssim \seminorm{\ul{a}_{T}}{\vThyb}^2 .
  \end{equation}
  \begin{equation} \label{eq:hho-rec-minus-cell}
    \seminorm{1,T}{\rec{T}{\K+1}\vThyb - v_T} \lesssim \HT^{-1/2} \norm{\partial T}{v_{\partial T} - \lproj{\partial T}{\K} v_T} = s_T(\vThyb, \vThyb)^{1/2}.
  \end{equation}
\end{lemma}
\begin{proof}
  The proof for \eqref{eq:hho-cell-coercivity} can be found in \cite[Lemma 3, (3.1a)]{badia2026geometricmultigridsolvershybrid}. For \eqref{eq:hho-rec-minus-cell}, we take $w = \rec{T}{\K+1}\vThyb - v_T \in \Pk{\K+1}(T)$ in \eqref{eq:hho-rec}, then insert $\lproj{\partial T}{\K}$ in the face term, which is allowed since $\bm{n}_F \cdot \nabla w|_F \in \IPk{\K}(F)$ by \eqref{eq:face-space}, to obtain 
  \[
  \seminorm{1,T}{w}^2 = (v_{\partial T} - \lproj{\partial T}{\K} v_T,\, \bm{n}_{\partial T} \cdot \nabla w)_{\partial T}.
  \]
  A Cauchy--Schwarz inequality, the trace inequality \eqref{eq:trace-inequality-continuous} and an inverse inequality applied to $\nabla w$ then yield \eqref{eq:hho-rec-minus-cell}.
\end{proof}

We also state the following norm equivalence result. Its proof in the case of same-order \ac{hho} is standard and can be found, e.g.\ in \cite{DiPietro-Droniou-2020-HHO-Book}. For the mixed-order method with the face spaces \eqref{eq:face-space} and the stabilisation \eqref{eq:hho-stab}, it can be found in \cite[Lemma 3]{badia2026geometricmultigridsolvershybrid}.
\begin{lemma}[\ac{hho} seminorm equivalence] \label{lem:hho-norm-equivalence}
  If $\TH$ satisfies Assumption \ref{ass:coarse-mesh-regularity}, then, for all $\vThyb \in \UHhyb(T)$,
  \[
    \normHHOloc{\vThyb}^2 \lesssim \seminorm{\ul{a}_{T}}{\vThyb}^2 \lesssim \normHHOloc{\vThyb}^2.
  \]
\end{lemma}

Summing over $T \in \TH$, and using the definition of the global norms, Lemma \ref{lem:hho-norm-equivalence} immediately gives the global equivalence $\normHHO{\vHhyb} \simeq \norm{\ul{a}_H}{\vHhyb}$ for all $\vHhyb \in \UHhyb$, which is the form used in the remainder of this section.

\begin{lemma}[Positive definiteness] \label{lem:hho-pd}
  If $\TH$ satisfies Assumption \ref{ass:coarse-mesh-regularity}, then the bilinear form $\ul{a}_H$ is positive definite on $\UHhyb$, and the condensed bilinear form $\aH$ is positive definite on $\MH$.
\end{lemma}
\begin{proof}
  The statement for $\aH$ is \cite[Lemma 4]{badia2026geometricmultigridsolvershybrid}. For $\ul{a}_H$, if $\ul{a}_H(\vHhyb, \vHhyb) = 0$ then, by Lemma \ref{lem:hho-norm-equivalence}, $v_T$ is constant and equal to $v_{\partial T}$ on every $T \in \TH$, and the argument of \cite[Lemma 4]{badia2026geometricmultigridsolvershybrid} applies. Cells sharing a face carry the same constant, the constant vanishes on cells with a boundary face, and every cell is joined to one such by a chain of faces since $\dom$ is connected. Hence $\vHhyb = 0$.
\end{proof}

We now define the HHO interpolation operator $\JThyb : H^1(T) \to \UHhyb(T)$ as follows:
\begin{equation}
  \label{eq:hho-interp-def}
  \JThyb v \doteq \bigl(\lproj{T}{\K+1} v,\, \lproj{\partial T}{\K} v\bigr).
\end{equation}
Its global counterpart $\JHhyb : H^1_0(\dom) \to \UHhyb$ is defined by applying $\JThyb$ element-wise. This is well defined since the trace of $v \in H^1_0(\dom)$ on an interior face is single-valued. Note that, for $v \in H^1_0(\dom)$, we have $\lproj{F}{\K} v = 0$ on every boundary face $F \in \FHb$, so that indeed $\JHhyb v \in \UHhyb$. We denote by $\trH$ the trace operator returning the skeleton component of a hybrid function, i.e.\ $\trH \vHhyb \doteq v_{\FH}$ and, locally, $\trH \vThyb \doteq v_{\partial T}$.

As a consequence of \eqref{eq:hho-intermediate-trace-bound-2}, \eqref{eq:l2-approx-trace} and \eqref{eq:l2-approx}, the HHO interpolation is bounded (see \cite[Lemma 8, (3.6a)]{badia2026geometricmultigridsolvershybrid}): for any $v \in H^1(T)$,
\begin{equation}
  \label{eq:hho-interpolation-boundedness}
  \normHHOloc{\JThyb v} \lesssim \seminorm{1,T}{v}.
\end{equation}
A fundamental property of the HHO reconstruction is that it commutes with the elliptic projector. For all $v \in H^1(T)$ we have $\rec{T}{\K+1}\bigl(\JThyb v\bigr) = \eproj{T}{\K+1} v$, where $\eproj{T}{\K+1} : H^1(T) \to \Pk{\K+1}(T)$ is the elliptic projector defined in \eqref{eq:elliptic-projector}. For the face spaces \eqref{eq:face-space}, the proof for the same-order \ac{hho} formulation can be found in \cite[Lemma 2]{Yemm2022ANA} and extends verbatim to the mixed-order method used here \cite{badia2026geometricmultigridsolvershybrid}. As observed in \cite[Remark 4]{Yemm2022ANA}, the proof relies on the fact that $\nabla \Pk{\K+1}(T) \cdot \bm{n}_F|_F \subset \IPk{\K}(F)$.

As a consequence, the approximation properties of $\rec{T}{\K+1} \circ \JThyb$ reduce to standard estimates for the elliptic projector on shape-regular polytopes, given by \eqref{eq:h1-approx} and \eqref{eq:h1-approx-trace}.

\subsection{HHO prolongation}

As our prolongation $\IHh = \IHhhho$ $: \MH \to \Uh$, we consider an operator inspired by the ones used in \cite{badia2026geometricmultigridsolvershybrid}. It first maps the skeleton variable to a \ac{dg} function on the coarse mesh, then applies the fine interpolant and the averaging operator $\Wh$ of Section \ref{sec:dg} to enforce conformity. Since there is more than one choice for the coarse reconstruction, we will formulate a general framework for the analysis of $\IHhhho$ like the one in \cite{badia2026geometricmultigridsolvershybrid}.

\begin{definition}
  \label{def:hho-prolongation}
  Let $\Theta_H: \MH \to \Pk{\K+1}(\TH)$ be a coarse reconstruction operator mapping skeleton functions to discontinuous coarse bulk polynomials, and let $\IHhdg = \Wh \Ih$ be the \ac{dg} prolongation \eqref{eq:prolongation-def}, which is defined on $\Pk{k}(\TH)$ for any $k \geq 1$. We define the HHO prolongation as
  \[
    \IHhhho \doteq \IHhdg \Theta_H = \Wh \Ih \Theta_H.
  \]
  In particular, we will consider the choices $\Theta_H = \iscH$ and $\Theta_H = \rec{H}{\K+1} \circ \iscHhyb$.
\end{definition}

We first collect the properties of the two reconstructions $\Theta_H$ that the analysis needs, starting with two results from \cite{badia2026geometricmultigridsolvershybrid} on the energy of bubble functions.

\begin{lemma} \label{lem:bubble-function-bound}
  Let $\ul{\delta}_H = (\delta_H, 0) \in \UHhyb$ be a bubble function on $\UHhyb$. Then, we have the following bound:
  \[
    \norm{L^2(\dom)}{\delta_H} \lesssim \H \norm{\ul{a}_H}{\ul{\delta}_H}.
  \]
\end{lemma}
\begin{proof}
  The proof can be found in \cite[Lemma 10, (3.10a)]{badia2026geometricmultigridsolvershybrid}.
\end{proof}

The next result is the $H^1$ counterpart of \cite[Lemma 10, (3.10b)]{badia2026geometricmultigridsolvershybrid}, which assumes $z \in H^2(\dom)$ and gains a factor $\H$. The orthogonality argument of the proof is the same, but the consistency term is bounded through the elliptic projector rather than through the consistency error.

\begin{lemma}
  \label{lem:continuous-bubble-bound}
  Let $z \in H^1_0(\Omega)$. Define two hybrid functions:
  \[
    \uHhyb \doteq \JHhyb z = (\lproj{\TH}{\K+1} z, \lproj{\FH}{\K} z)\ , \quad \vHhyb \doteq \iscHhyb  \trH  \JHhyb z = (\iscH \lproj{\FH}{\K} z, \lproj{\FH}{\K} z).
  \]
  Then, we have the following bound:
  \[
    \norm{\ul{a}_H}{\uHhyb - \vHhyb} \lesssim \seminorm{1,\dom}{z}.
  \]
\end{lemma}
\begin{proof}
  We define the bubble function $\ul{\delta}_H = \uHhyb - \vHhyb = (\delta_H, 0)$. Then, by a-orthogonality of $\iscH$ with bubble functions, we have
  \[
    \norm{\ul{a}_H}{\ul{\delta}_H}^2 = \ul{a}_H(\ul{\delta}_H, \ul{\delta}_H) = \ul{a}_H(\uHhyb, \ul{\delta}_H) = \left[\sum_{T \in \TH} \left( \nabla \rec{T}{\K+1} \uHhyb, \nabla \rec{T}{\K+1} \ul{\delta}_H \right)_T\right] + \sum_{T \in \TH} s_T(\uHhyb, \ul{\delta}_H).
  \]
  The consistency term can be bounded by remembering that $\rec{T}{\K+1} \uHhyb = \rec{T}{\K+1} \JHhyb z = \eproj{T}{\K+1} z$ and thus
  \begin{align*}
    \sum_{T \in \TH} \left( \nabla \rec{T}{\K+1} \uHhyb, \nabla \rec{T}{\K+1} \ul{\delta}_H \right)_T
    &= \sum_{T \in \TH} \left( \nabla \eproj{T}{\K+1} z, \nabla \rec{T}{\K+1} \ul{\delta}_H \right)_T
    = \sum_{T \in \TH} \left( \nabla z, \nabla \rec{T}{\K+1} \ul{\delta}_H \right)_T \\
    &\leq \seminorm{1,\dom}{z} \seminorm{1,\TH}{\rec{H}{\K+1} \ul{\delta}_H}
    \leq \seminorm{1,\dom}{z} \norm{\ul{a}_H}{\ul{\delta}_H}
  \end{align*}
  where we have used that $(\nabla \eproj{T}{\K+1} z, \nabla w)_T = (\nabla z, \nabla w)_T$ for all $w \in \Pk{\K+1}(T)$, then a Cauchy-Schwarz inequality.
  The stabilisation term can be bounded as follows:
  \begin{equation*}
    \sum_{T \in \TH} s_T(\uHhyb, \ul{\delta}_H)
    \leq \left(\sum_{T \in \TH} s_T(\JHhyb z, \JHhyb z)\right)^{1/2} \left(\sum_{T \in \TH} s_T(\ul{\delta}_H, \ul{\delta}_H)\right)^{1/2}
    \lesssim \seminorm{1,\dom}{z} \norm{\ul{a}_H}{\ul{\delta}_H}
  \end{equation*}
  where we have first used a Cauchy-Schwarz inequality, and then $s_T(\JThyb z, \JThyb z) = \HT^{-1} \norm{\partial T}{\lproj{\partial T}{\K}(\lproj{T}{\K+1}z - z)}^2 \lesssim \HT^{-1} \norm{\partial T}{z - \lproj{T}{\K+1}z}^2 \lesssim \seminorm{1,T}{z}^2$, by $L^2$-stability of $\lproj{\partial T}{\K}$ and \eqref{eq:l2-approx-trace} with $s = 1$.
  Putting everything together and dividing by $\norm{\ul{a}_H}{\ul{\delta}_H}$ (the case $\ul{\delta}_H = 0$ being trivial) concludes the proof.
\end{proof}

The following approximation result is the $H^1$ counterpart of \cite[Lemmas 15 and 17]{badia2026geometricmultigridsolvershybrid}, which are stated one order higher for $H^2$ functions.

\begin{lemma} \label{lem:hho-approximation}
  Given $\vh \in \Uh$, let $\lH \doteq \lproj{\FH}{\K} \vh \in \MH$. Then, for $r = 0, 1$ and both choices $\Theta_H = \iscH$ and $\Theta_H = \rec{H}{\K+1} \circ \iscHhyb$,
  \begin{equation}
    \label{eq:hho-approx}
    \norm{H^r(\TH)}{\vh - \Theta_H \lH} \lesssim \H^{1-r} \seminorm{1,\dom}{\vh}.
  \end{equation}
\end{lemma}
\begin{proof}
  Using a triangle inequality, we can write, for $r = 0, 1$,
  \[
    \norm{H^r(\TH)}{\vh - \Theta_H \lH} \leq \norm{H^r(\TH)}{\vh - \lproj{\TH}{\K+1} \vh} + \norm{H^r(\TH)}{\lproj{\TH}{\K+1} \vh - \Theta_H \lH}.
  \]
  By the approximation properties of the $L^2$-orthogonal projector in Equation \eqref{eq:l2-approx}, applied element-wise on $\TH$ with $s = 1$ and $m = r$, we have that
  \[
    \norm{H^r(\TH)}{\vh - \lproj{\TH}{\K+1} \vh} \lesssim \H^{1-r} \seminorm{1,\dom}{\vh}.
  \]
  We now tackle the second term. First consider the case $\Theta_H = \iscH$, and define the two hybrid functions of Lemma \ref{lem:continuous-bubble-bound} with $z = \vh \in H^1_0(\dom)$ as $\uHhyb = \JHhyb \vh$ and $\vHhyb = \iscHhyb  \trH  \JHhyb \vh$. The bubble function $\ul{\delta}_H = \uHhyb - \vHhyb$ satisfies that $\delta_H = \lproj{\TH}{\K+1} \vh - \iscH \lproj{\FH}{\K} \vh = \lproj{\TH}{\K+1} \vh - \iscH \lH$. Then, using Lemmas \ref{lem:bubble-function-bound} and \ref{lem:continuous-bubble-bound}, we have that
  \[
    \norm{L^2(\dom)}{\delta_H} \lesssim \H \norm{\ul{a}_H}{\ul{\delta}_H} \lesssim \H \seminorm{1,\dom}{\vh}.
  \]
  Moreover, by Lemma \ref{lem:hho-norm-equivalence} for the seminorm part and the previous bound with $\H \leq 1$ for the $L^2$ part, together with Lemma \ref{lem:continuous-bubble-bound},
  \[
    \norm{H^1(\TH)}{\delta_H} \lesssim \norm{\ul{a}_H}{\ul{\delta}_H} \lesssim \seminorm{1,\dom}{\vh}.
  \]
  This concludes the proof for the case $\Theta_H = \iscH$. For the case $\Theta_H = \rec{H}{\K+1} \circ \iscHhyb$, we start by using a triangle inequality to write
  \[
    \norm{H^r(\TH)}{\vh - \rec{H}{\K+1} \circ \iscHhyb \lH} 
    \leq \norm{H^r(\TH)}{\vh - \iscH \lH} + \norm{H^r(\TH)}{\iscH \lH - \rec{H}{\K+1} \circ \iscHhyb \lH}.
  \]
  We can bound the first term using the previous case. For the second term, we apply \eqref{eq:hho-rec-minus-cell} to $\vHhyb = \iscHhyb \lH$, for which $v_T = \isc{T}\lH$ and $v_{\partial T} = \lH|_{\partial T}$. Recalling the definition \eqref{eq:hho-stab} of the stabilisation, and summing over $T \in \TH$, this gives
  \[
    \seminorm{1,\TH}{\iscH \lH - \rec{H}{\K+1} \circ \iscHhyb \lH}^2
    \lesssim \sum_{T \in \TH} s_T(\iscHhyb \lH, \iscHhyb \lH)
    \leq \norm{\ul{a}_H}{\iscHhyb \lH}^2.
  \]
  Then, with the same $\uHhyb$, $\vHhyb$ and $\ul{\delta}_H$ defined above, we have $\iscHhyb \lH = \vHhyb = \uHhyb - \ul{\delta}_H$, so that a triangle inequality gives
  \[
    \norm{\ul{a}_H}{\iscHhyb \lH} \leq \norm{\ul{a}_H}{\ul{\delta}_H} + \norm{\ul{a}_H}{\uHhyb} \lesssim \seminorm{1,\dom}{\vh}, 
  \]
  where the bubble term is bounded by Lemma \ref{lem:continuous-bubble-bound}, and the term $\uHhyb = \JHhyb \vh$ by the continuity of $\ul{a}_H$ given by Lemma \ref{lem:hho-norm-equivalence} together with the boundedness \eqref{eq:hho-interpolation-boundedness}.
  Finally, we note that, by definition of $\rec{H}{\K+1}$, $(\iscH \lH - \rec{H}{\K+1} \circ \iscHhyb \lH, 1)_T = 0$ for all $T \in \TH$, and thus we can apply a Poincaré inequality \eqref{eq:poincare-inequality} to obtain
  \[
    \norm{L^2(\dom)}{\iscH \lH - \rec{H}{\K+1} \circ \iscHhyb \lH} \lesssim \H \seminorm{1,\TH}{\iscH \lH - \rec{H}{\K+1} \circ \iscHhyb \lH}.
  \]
  Indeed, for $r = 1$ the full $H^1(\TH)$ norm is controlled by the seminorm bound together with the $L^2$ bound and $\H \leq 1$, while for $r = 0$ the $L^2$ bound alone gives the claim.
\end{proof}

\begin{lemma}[Control of the jumps] \label{lem:hho-jumps}
  For all $\vHhyb = (v_{\TH}, v_{\FH}) \in \UHhyb$,
  \[
    \sum_{F \in \FH} H_F^{-1} \norm{F}{\jump{v_{\TH}}}^2 + \sum_{F \in \FH} H_F^{-1} \norm{F}{\jump{\rec{H}{\K+1} \vHhyb}}^2 \lesssim \normHHO{\vHhyb}^2.
  \]
\end{lemma}
\begin{proof}
  The proof can be found in \cite[Lemma 14]{badia2026geometricmultigridsolvershybrid}, where the bound is stated with the uniform coarse mesh size $\H$ in place of $H_F$. The two are equivalent by Assumption~\ref{ass:coarse-mesh-regularity}(i).
\end{proof}

\begin{lemma} \label{lem:hho-theta-properties}
  For both choices of $\Theta_H: \MH \to \Pk{\K+1}(\TH)$ in Definition \ref{def:hho-prolongation}, namely $\Theta_H = \iscH$ and $\Theta_H = \rec{H}{\K+1} \circ \iscHhyb$, we have that, for all $\lH \in \MH$, the following bounds hold:
  \begin{align}
    \seminorm{1,\TH}{\Theta_H \lH}^2 &\lesssim \norm{\aH}{\lH}^2, \label{eq:hho-theta-h1} \\
    \sum_{F \in \FH} H_F^{-1} \norm{F}{\jump{\Theta_H \lH}}^2 &\lesssim \norm{\aH}{\lH}^2. \label{eq:hho-theta-jumps}
  \end{align}
\end{lemma}
\begin{proof}
  We first prove \eqref{eq:hho-theta-h1}. For $\Theta_H = \rec{H}{\K+1} \circ \iscHhyb$, the gradient of the reconstruction is the first term of the local bilinear form \eqref{eq:hho-form-local}, so that $\seminorm{1,T}{\rec{T}{\K+1} \iscHhyb \lH}^2 \leq a_T(\lH, \lH)$, and summing over $T \in \TH$ gives the bound. For $\Theta_H = \iscH$, the bound follows from \eqref{eq:hho-cell-coercivity} applied to $\iscHhyb \lH$ and summed over $T \in \TH$. 
  Then Lemma \ref{lem:hho-jumps} applied to $\vHhyb = \iscHhyb \lH$, together with the norm equivalence of Lemma \ref{lem:hho-norm-equivalence}, gives \eqref{eq:hho-theta-jumps} for both choices of $\Theta_H$.
\end{proof}

The following three lemmas are the counterparts of Lemmas \ref{lem:oswald-dg} and \ref{lem:oswald-prolongation}. The reconstruction degree $\K+1$ may exceed the fine polynomial degree $\k$, and does so whenever $\K \geq \k$, which includes the choice $\K = \k$ made in our experiments. This case is covered by the bounds of Lemmas \ref{lem:oswald}--\ref{lem:oswald-prolongation}, which hold for an arbitrary coarse polynomial degree $k$.

\begin{lemma} \label{lem:hho-h1}
  For both choices of $\Theta_H$ in Definition \ref{def:hho-prolongation} and all $\lH \in \MH$, the prolongation operator $\IHhhho$ satisfies
  \begin{equation}
    \norm{H^1(\dom)}{\IHhhho \lH}^2 \lesssim \frac{\H}{\h} \norm{\aH}{\lH}^2.
  \end{equation}
\end{lemma}
\begin{proof}
  Set $\vH \doteq \Theta_H \lH \in \Pk{\K+1}(\TH)$, so that $\IHhhho \lH = \IHhdg \vH$. Since $\IHhdg \vH \in H^1_0(\dom)$, a Poincar\'e inequality and the triangle inequality give
  \[
    \norm{H^1(\dom)}{\IHhhho \lH} \lesssim \seminorm{1,\dom}{\IHhdg \vH} \leq \seminorm{1,\TH}{\vH - \IHhdg \vH} + \seminorm{1,\TH}{\vH}.
  \]
  By Lemma \ref{lem:oswald} with $r = 1$ and $k = \K + 1$, together with $\h \leq H_T$,
  \[
    \norm{H^1(\dom)}{\IHhhho \lH}^2 \lesssim \sum_{F \in \FH} \h^{-1} \norm{F}{\jump{\Theta_H \lH}}^2 + \seminorm{1,\TH}{\Theta_H \lH}^2.
  \]
  Since $H_F \simeq \H$, the jump term satisfies
  \[
    \sum_{F \in \FH} \h^{-1} \norm{F}{\jump{\Theta_H \lH}}^2
    = \sum_{F \in \FH} \frac{H_F}{\h}\, H_F^{-1} \norm{F}{\jump{\Theta_H \lH}}^2
    \lesssim \frac{\H}{\h} \sum_{F \in \FH} H_F^{-1} \norm{F}{\jump{\Theta_H \lH}}^2
    \lesssim \frac{\H}{\h}\, \norm{\aH}{\lH}^2,
  \]
  by \eqref{eq:hho-theta-jumps}, and the volume term is bounded by \eqref{eq:hho-theta-h1}.
\end{proof}

\begin{lemma} \label{lem:hho-prolongation}
  Given $\vh \in \Uh$, let $\lH \doteq \lproj{\FH}{\K} \vh \in \MH$. Then, for both choices of $\Theta_H$ in Definition \ref{def:hho-prolongation}, the prolongation operator $\IHhhho$ satisfies
  \begin{equation}
    \norm{L^2(\dom)}{\vh - \IHhhho\lH}^2 \lesssim \H^2 \seminorm{1,\dom}{\vh}^2 \quad \text{and} \quad
    \norm{H^1(\dom)}{\vh - \IHhhho\lH}^2 \lesssim \frac{\H}{\h} \seminorm{1,\dom}{\vh}^2.
  \end{equation}
\end{lemma}
\begin{proof}
  Since $\IHhhho \lH = \IHhdg (\Theta_H \lH)$, this is Lemma \ref{lem:oswald-prolongation} applied to $\vH = \Theta_H \lH$, whose hypothesis \eqref{eq:coarse-approx} is Lemma \ref{lem:hho-approximation}.
\end{proof}

\begin{lemma} \label{lem:hho-bound}
  Given $\vh \in \Uh$, let $\lH = \lproj{\FH}{\K} \vh \in \MH$. Then we have that
  \[
    \norm{\aH}{\lH} \lesssim \seminorm{1,\dom}{\vh}.
  \]
\end{lemma}
\begin{proof}
  Consider the hybrid function $\JHhyb \vh = (\lproj{\TH}{\K+1} \vh, \lproj{\FH}{\K} \vh) \in \UHhyb$, whose skeleton component is precisely $\lH$. By the energy minimisation property \eqref{eq:energy-minimisation}, the continuity of $\ul{a}_H$ given by Lemma \ref{lem:hho-norm-equivalence} and the boundedness \eqref{eq:hho-interpolation-boundedness} of the HHO interpolator, we obtain
  \[
    \norm{\aH}{\lH}^2 = \aH(\lH, \lH) \leq \ul{a}_H(\JHhyb \vh, \JHhyb \vh)
    \lesssim \normHHO{\JHhyb \vh}^2 \lesssim \seminorm{1,\dom}{\vh}^2,
  \]
  which concludes the proof.
\end{proof}
\subsection{Main result}

\begin{theorem}
  \label{th:hho-main}
  The \ac{hho} discretisation, with the smoother described in Section \ref{sec:smoother} and the prolongation $\IHh = \IHhhho$ of Definition \ref{def:hho-prolongation}, fulfils the assumptions in Lemma \ref{lem:auxspace-concrete}.
\end{theorem}
\begin{proof}
  That the smoother satisfies \eqref{eq:auxspace-abstract-smoother} is given by Corollary \ref{cor:point-jacobi} or \ref{cor:vertex-patch}, depending on the smoother at hand.

  Let $\lH \in \MH$. By Lemma \ref{lem:hho-h1},
  \[
    \norm{\ah}{\IHhhho \lH}^2 = \seminorm{1,\dom}{\IHhhho \lH}^2 \lesssim \frac{\H}{\h} \norm{\aH}{\lH}^2,
  \]
  which proves \eqref{eq:auxspace-abstract-prolongation}.

  We now turn to \eqref{eq:auxspace-abstract-decomposition}. Given $\vh \in \Uh$, we choose as coarse approximant $\lH \doteq \lproj{\FH}{\K} \vh \in \MH$. Lemma \ref{lem:hho-prolongation}, together with $\seminorm{1,\dom}{\vh} = \norm{\ah}{\vh}$, gives
  \[
    \norm{L^2(\dom)}{\vh - \IHhhho \lH}^2 \lesssim \H^2 \norm{\ah}{\vh}^2 \quad \text{ and } \quad \norm{H^1(\dom)}{\vh - \IHhhho \lH}^2 \lesssim \frac{\H}{\h} \norm{\ah}{\vh}^2.
  \]
  Finally, Lemma \ref{lem:hho-bound} gives $\norm{\aH}{\lH} \lesssim \seminorm{1,\dom}{\vh} = \norm{\ah}{\vh}$, which proves \eqp{eq:auxspace-abstract-decomposition}{2,3} and \eqp{eq:auxspace-abstract-decomposition}{1}, respectively.
\end{proof}

\begin{remark} \label{rem:hho-final}
  Analogously to the \ac{dg} case, Theorem \ref{th:hho-main} together with Lemma \ref{lem:auxspace-concrete} yields the final bound
  \[
    \kappa(P_{\mathrm{add,HHO}}^{-1} \Ah) \lesssim \Bigl( 1 + \lambda_1 \frac{\H}{\h} \Bigr)\Bigl(\frac{\H^2}{\h^2} + \frac{1}{\lambda_0}\Bigr),
  \]
  for the point-Jacobi and vertex-star patch smoothers considered in this work, for which $c_{D,0} \simeq \h^{-2}$, $c_{D,1} = 0$ and $c_D^2 \simeq 1$, and again $O(\H^3/\h^3)$ for an exact coarse solve. 
\end{remark}

\section{Coarse solvers and geometric multigrid}
\label{sec:multilevel}

Sections \ref{sec:smoother} to \ref{sec:hho} supply all ingredients of Lemma
\ref{lem:auxspace-concrete} except the ones related to the accuracy of the coarse solver, that is $\{\lambda_0, \lambda_1\}$ in Assumption \ref{ass:coarse-solver}. 

A first approach is to solve the coarse problem exactly, which gives $\BH = \AH^{-1}$ and thus $\lambda_0 = \lambda_1 = 1$. This requires the exact solution of a problem that is smaller than the original one by a factor $(\H/\h)^{d}$. In real applications, this is often too costly, motivating the use of iterative solvers for the coarse problem.
Therefore, scalable iterative solvers for polytopal methods on agglomerated meshes are required to efficiently solve the coarse problem. 

In this work, we will employ \ac{gmg} solvers for this purpose. \ac{gmg} is a family of multilevel solvers defined on hierarchies of nested meshes, where each level corresponds to a different resolution of the computational domain. Recent advances have demonstrated the effectiveness of \ac{gmg} on agglomeration-based hierarchies for both \ac{hho} \cite{badia2026geometricmultigridsolvershybrid} and \ac{dg} \cite{Antonietti2017Multigrid} methods. These works prove bounds for the spectrum of $\BH \AH$ when $\BH$ is the multigrid V-cycle with appropriate prolongation operators and smoothers. These bounds are independent of the number of levels of the hierarchy and the size $\H$ of the finest problem (here our coarse problem $\AH$). Instead, these bounds mostly depend on the coarsening ratio, which can be kept constant throughout the hierarchy.

\subsection{Geometric multigrid coarse solver} \label{sec:gmg-solver}

We now briefly discuss the multigrid algorithm, including the choice of prolongation operators and smoothers for the \ac{gmg} coarse solver that will be used in the numerical experiments. We refer to \cite{badia2026geometricmultigridsolvershybrid,Antonietti2017Multigrid} for details.

The mesh $\TH$, which is the coarse mesh of the two-level method, is the finest mesh seen by the coarse solver. Starting from it, we construct by further agglomeration a hierarchy of nested meshes $\{\mathcal{T}_\ell\}_{\ell=1}^L$, with $\mathcal{T}_1 = \TH$ at the finest level $\ell = 1$ and $\mathcal{T}_L$ at the coarsest level $\ell = L$. For either \ac{dg} or \ac{hho}, we construct on each level $\ell$ a discrete space $V_\ell$ and its associated operator $A_\ell$, analogously to the constructions in Sections \ref{sec:dg} and \ref{sec:hho}, respectively. 

For $\ell < L$, let $I_\ell : V_{\ell+1} \to V_\ell$ denote the prolongation operator from level $\ell+1$ to level $\ell$, and $I_\ell^{*}$ its adjoint. For \ac{dg}, we follow \cite{Antonietti2017Multigrid} and take $I_\ell$ to be the injection operator, since we have $V_{\ell+1} \subset V_\ell$. For \ac{hho}, more complicated prolongation operators are required due to the hybrid nature of the method. In our experiments, we consider $I_\ell$ as described in \cite[Section 4.2]{badia2026geometricmultigridsolvershybrid}, which given a hybrid coarse function $\mu$, constructs its fine-level prolongation by averaging the cell values of its inverse static condensation $\mathcal{U}\mu$ over fine faces. The construction is analogous to the one in Definition \ref{def:hho-prolongation}, with the averaging operator $\Wh$ replaced by an $L^2$ face-averaging operator. Operators of this kind are common in the multilevel literature for hybrid methods, and can be found in, e.g.\ \cite{Di_Pietro-Dong-Kanschat-Rupp-2024-Multigrid,DiPietro2021,Lu2022}. 

We turn to the choice of smoothers. In \cite{Antonietti2017Multigrid}, a standard point-Jacobi smoother is used. In this work, we will instead employ the block-Jacobi smoothers described in \cite[Section 4.3]{badia2026geometricmultigridsolvershybrid} for both the \ac{hho} and \ac{dg} implementations. These are defined as follows: on a level $\ell < L$, let $D_\ell$ be defined as in \eqref{eq:smoother-def-inverse}, associated with the set $\mathcal{D}_\ell$ of macro-vertices of the coarse mesh, that is, vertices shared by multiple agglomerated cells (or interfaces in the sense of Section \ref{sec:auxspace}). The subdomains then correspond to the vertex stars around each macro-vertex. From $D_\ell$ we can define its associated damped Richardson smoother on level $\ell$, namely 
\begin{equation}
  \label{eq:richardson-smoother}
  S_{\ell,n_s,\omega} \doteq (\identity - \omega D_\ell^{-1} A_\ell)^{n_s},
\end{equation}
where $\identity$ is the identity on $V_\ell$, $\omega \in (0,1)$ is the damping parameter, and $n_s$ is the number of smoothing steps. 

With these two ingredients, we can define the \ac{gmg} V-cycle through its error propagation operator. Let $B_L = A_L^{-1}$ be the exact solver on the coarsest level. Then for the levels $\ell = L-1, \dots, 1$, we have
\begin{equation}
  \label{eq:gmg-v-cycle}
  \identity - B_{\ell} A_\ell \doteq S_{\ell,n_s,\omega} \, (\identity - I_\ell B_{\ell+1} I_\ell^{*} A_\ell) \, S_{\ell,n_s,\omega}.
\end{equation}
We can then set $\BH \doteq B_1$. To describe the two multilevel algorithms below in a single notation, we add the conforming problem of \sect{sec:preliminaries} as level $\ell = 0$ of the hierarchy, setting $\mathcal{T}_0 \doteq \Th$, $V_0 \doteq \Uh$, $A_0 \doteq \Ah$, $D_0 \doteq \Dh$ and $I_0 \doteq \IHh$, so that the full hierarchy is $\{\mathcal{T}_\ell\}_{\ell=0}^L$ and the \ac{gmg} coarse solver acts on levels $1$ to $L$. This is the numbering used in \sect{sec:results}.

\subsection{The additive algorithm} \label{sec:additive}

With the \ac{gmg} V-cycle as coarse solver, the additive preconditioner \eqref{eq:auxspace-preconditioner} reads
\[
  P_{\mathrm{add}}^{-1} = I_0 B_1 I_0^{*} + D_0^{-1},
\]
an $(L+1)$-level algorithm in which level $0$ is treated additively and levels $1$ to $L$ multiplicatively, by the V-cycle. Its condition number follows from Lemma \ref{lem:auxspace-concrete} and the contraction factor of the V-cycle. For a convergent symmetric multigrid iteration, if its error propagation operator is $a_H$-self-adjoint and positive semi-definite with contraction factor $\rho < 1$, then
\[
  (1-\rho)\,(\AH^{-1}\mu_H, \mu_H) \;\leq\; (\BH \mu_H, \mu_H) \;\leq\; (\AH^{-1}\mu_H, \mu_H) \quad \forall \mu_H \in \UH,
\]
that is, $\lambda_1 \leq 1$ and $\lambda_0 = 1-\rho$, so that
\[
  \kappa(P_{\mathrm{add}}^{-1}\Ah) \lesssim \Bigl(c_D^2 + \frac{\H}{\h}\Bigr)\Bigl(\frac{\H^2}{\h^2} + \frac{1}{1-\rho}\Bigr).
\]

\subsection{A fully multiplicative algorithm} \label{sec:multiplicative}

The same ingredients give a fully multiplicative variant by treating level $0$ as one more level of the V-cycle. Let $S_{0,n_s,\omega} \doteq (\identity - \omega D_0^{-1} A_0)^{n_s}$ be the damped Richardson smoother \eqref{eq:richardson-smoother} on level $0$, which we also denote $S_{h,n_s,\omega}$ in \sect{sec:results}. Carrying the recursion \eqref{eq:gmg-v-cycle} one level further defines the multiplicative auxiliary space preconditioner $P_{\mathrm{mult}}^{-1} \doteq B_0$, that is,
\begin{equation}
  \label{eq:auxspace-preconditioner-multiplicative}
  \identity - P_{\mathrm{mult}}^{-1} A_0 \doteq S_{0,n_s,\omega} \, (\identity - I_0 B_1 I_0^{*} A_0) \, S_{0,n_s,\omega}.
\end{equation}
Both algorithms therefore share the hierarchy, the transfer operators, the smoothers and the coarse solver of \sect{sec:gmg-solver}, and differ only in whether level $0$ is treated additively or multiplicatively. Although we do not provide the analysis for $P_{\mathrm{mult}}$, we will benchmark its performance alongside $P_{\mathrm{add}}$ in the numerical experiments.

\section{Numerical results}
\label{sec:results}

We now test the preconditioners numerically in two and three dimensions, on several mesh families and polynomial degrees, and report iteration counts as a function of the mesh size and of the coarsening ratio.
To this end, we introduce three families of two-dimensional mesh hierarchies.

\emph{Cartesian meshes.}
First, we consider mesh hierarchies for which the agglomeration is regular and the coarse cells have flat interfaces. To this end, we consider $d$-dimensional Cartesian meshes on the unit cuboid $\Omega = [0,1]^d$. The hierarchy is then built recursively agglomerating blocks of $2^d$ neighbouring cells. After $\ell$ agglomeration steps, the coarse mesh $\mathcal{T}_\ell$ is composed of general polytopes with $2d$ flat interfaces, each composed of $(\H/\h)^{d-1}$ flat $(d-1)$-dimensional cuboids, and the coarsening ratio is exactly $\H/\h = 2^{\ell}$. 

\emph{Simplexified Cartesian meshes with arbitrary agglomeration.}
Starting from the same Cartesian meshes, which we simplexify, we now let the agglomeration be arbitrary. On each level, connected clusters of cells are agglomerated into polygonal elements, which are not required to be convex or shape-regular. Finding optimal agglomeration strategies for polytopal meshes remains an open research question \cite{FEDER2025,ANTONIETTI2026}. Within this work, the colouring of the mesh into clusters is performed using a hierarchical clustering algorithm with complete-linkage clustering and optimal branch ordering \cite{barjoseph2001}, provided by the Julia package \texttt{Clustering.jl} \cite{Clustering.jl}. The hierarchy is then built by applying the clustering algorithm recursively. The target reduction in the number of cells is specified per direction by a factor $\theta \in (0,1)$, so that each level reduces the target number of cells by $\theta^{d}$ and thus the coarsening ratio between levels is $\H/\h \simeq \theta^{-\ell}$. Due to the nature of these algorithms, the actual cluster sizes vary across the different agglomerates, but are on average close to the target. We consider $\theta = 0.5$, for which the measured (mean) ratios are within reasonable proximity of the nominal $2^{\ell}$.

\emph{Complicated geometries with arbitrary agglomeration.}
The third family starts from meshes $\Th$ of Australia, generated using the Julia package \texttt{GeoGmsh.jl}, which extracts the contour from publicly available shapefiles and then meshes the resulting geometry using \texttt{Gmsh} \cite{gmsh}. The coarsening is performed in the same way as in the previous example, but now using a k-means clustering algorithm provided by the Julia package \texttt{Clustering.jl} \cite{Clustering.jl}. For this set of experiments, we consider $\theta = 0.25$, so that $\H/\h \simeq 4$.

In three dimensions, we consider the three-dimensional counterparts of the Cartesian and simplexified meshes, with the same agglomeration strategies, on the unit cube $\Omega = [0,1]^3$.

\begin{figure}
  \centering
  \includegraphics[width=0.23\textwidth]{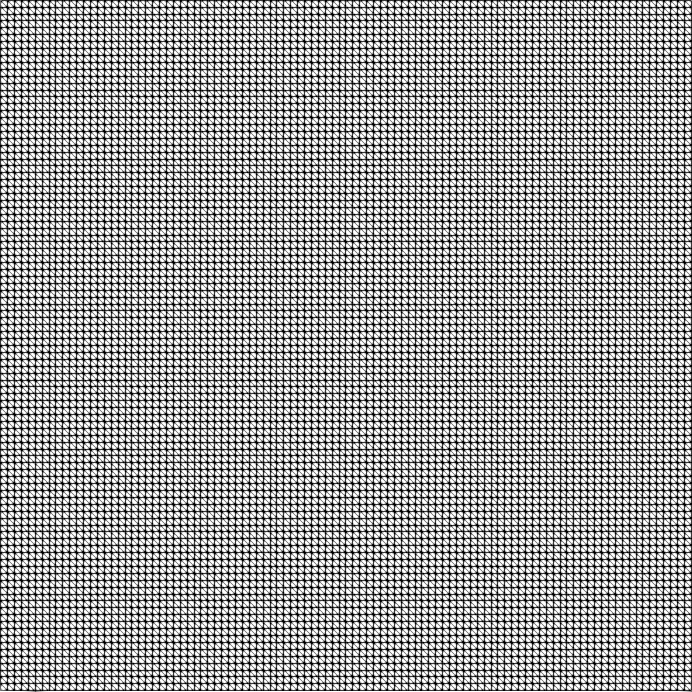}
  \includegraphics[width=0.23\textwidth]{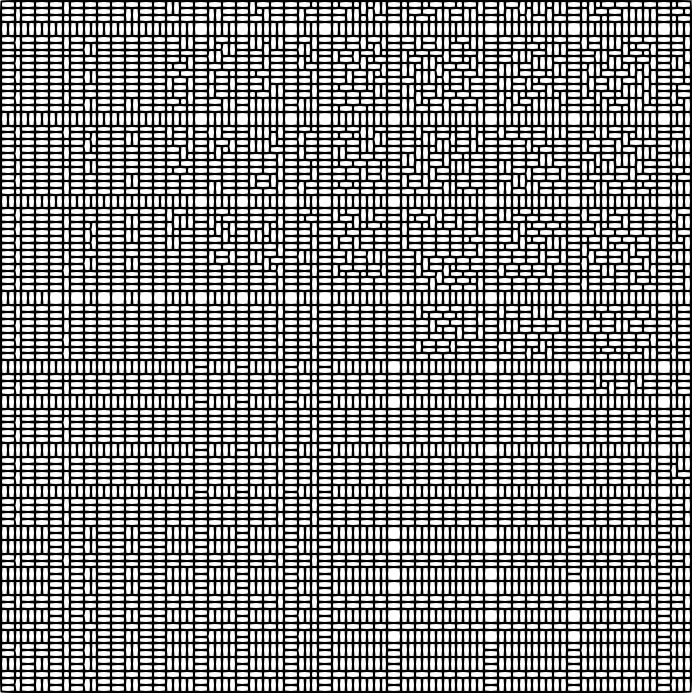}
  \includegraphics[width=0.23\textwidth]{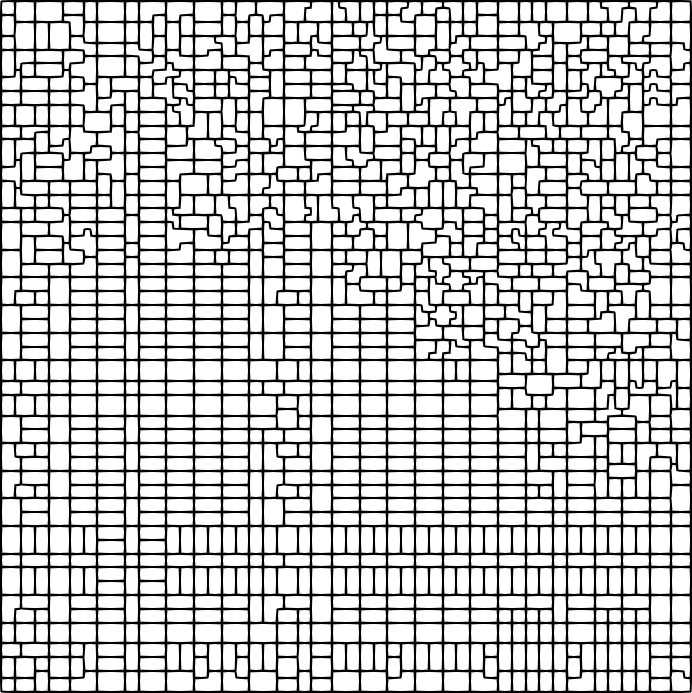}
  \includegraphics[width=0.23\textwidth]{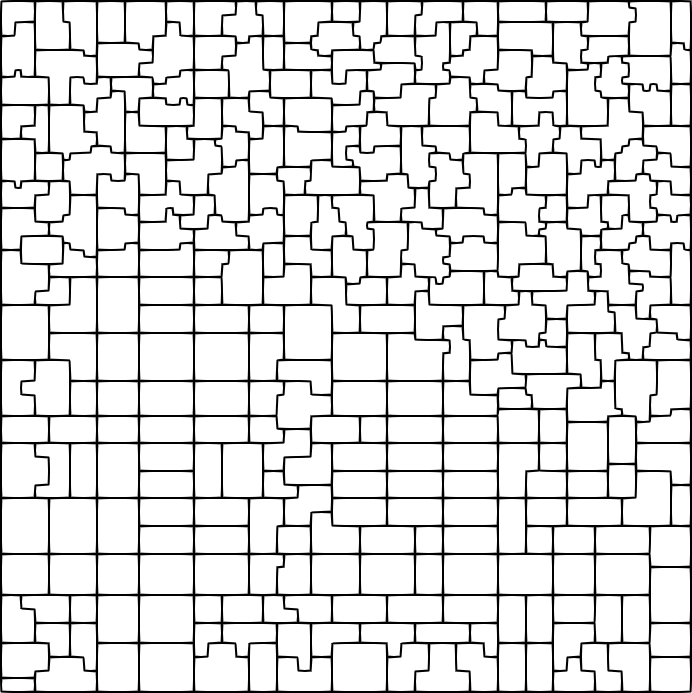} \\[2pt]
  \includegraphics[width=0.23\textwidth]{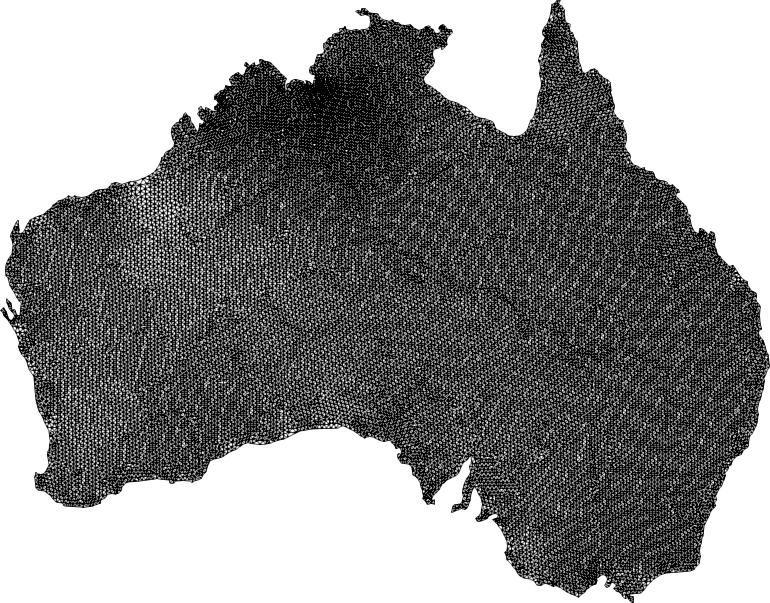}
  \includegraphics[width=0.23\textwidth]{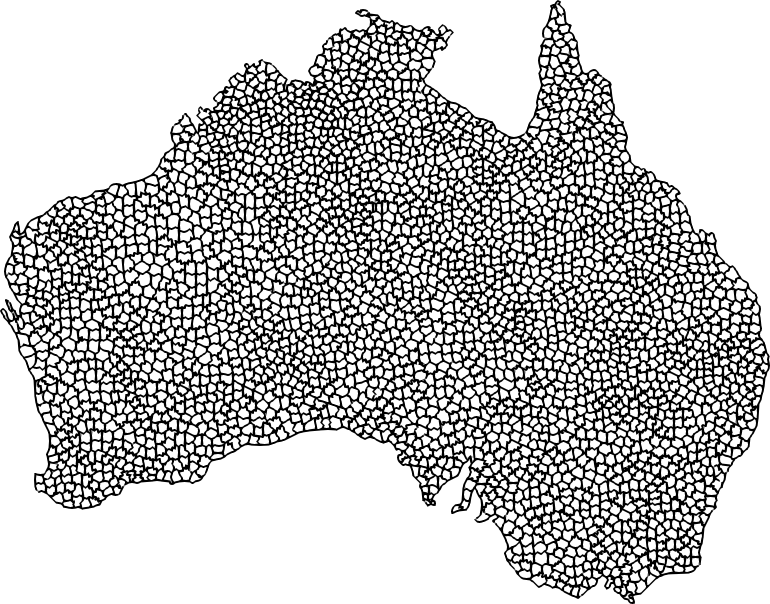}
  \includegraphics[width=0.23\textwidth]{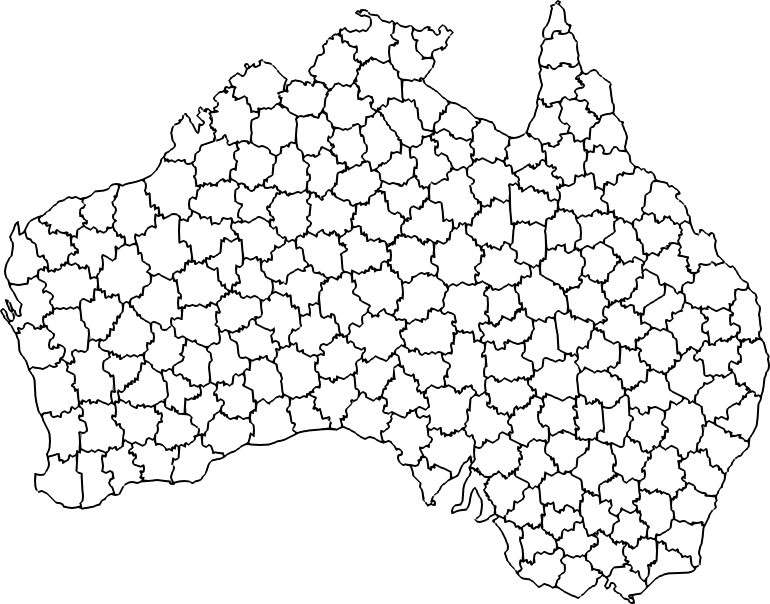}
  \includegraphics[width=0.23\textwidth]{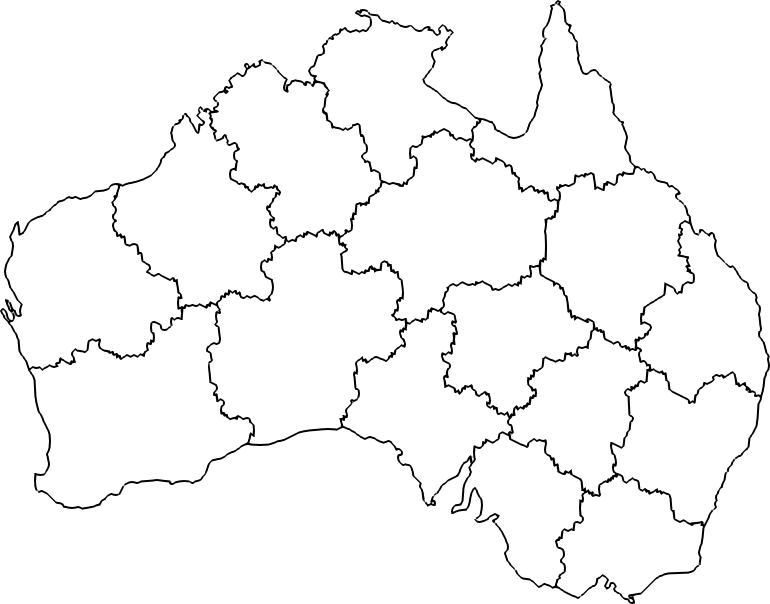}
  \caption{Agglomerated hierarchies for the \emph{Mid} simplexified Cartesian mesh with $\theta = 0.5$ (top) and the \emph{Mid} Australian mesh with $\theta = 0.25$ (bottom). Levels $\ell = 0$ to $3$, from left to right.}
  \label{fig:meshes}
\end{figure}

For each family of meshes, we generate several initial meshes of increasing resolution. For each initial mesh $\Th$, we build by agglomeration an associated hierarchy $\{\mathcal{T}_\ell\}_{\ell=0}^{L}$ with $\mathcal{T}_0 = \Th$ and $\mathcal{T}_1 = \TH$, so that $\{\mathcal{T}_\ell\}_{\ell=1}^{L}$ is the \ac{gmg} hierarchy of Section \ref{sec:multilevel}. Level $\ell = 0$ carries the conforming discretisation, and levels $\ell \geq 1$ the non-conforming ones. The number of cells at every agglomeration level, for each family and each initial mesh, is collected in Table \ref{tab:meshes}. Examples of these hierarchies are shown in Figure \ref{fig:meshes} for the simplexified Cartesian and Australian hierarchies.

On each initial mesh $\Th = \mathcal{T}_0$, we consider the conforming discretisation of Section \ref{sec:preliminaries} for a given polynomial degree $\k$. Then, on each level $\ell \geq 1$ of the hierarchy we build a space $V_\ell$ using either the \ac{dg} or \ac{hho} methods, as described in Sections \ref{sec:dg} and \ref{sec:hho}, respectively. In the tables these two variants are labelled CG-DG and CG-HHO. In both cases, the averaging operator $\Wh$ in the prolongation uses the uniform weights $\alpha_t(\nu) = |\mathcal{N}_h(\nu)|^{-1}$ of Remark \ref{rem:weights}, which is the classical Oswald interpolation operator. For the \ac{dg} method, we use the penalty parameter $\sigma = 10 (\K+1)^2$ throughout, and modify the \ac{hho} stabilisation weights to be $\simeq (\K+1)^2/\H$. These are typical choices in the literature to obtain a suitable scaling of the stabilisation term with respect to the polynomial order. Neither of these changes has an effect on the analysis, since polynomial order dependence was not tracked. 

In the experiments below, we report the number of preconditioned \ac{fgmres} iterations required to reduce the relative residual below $10^{-8}$. Using \ac{fgmres} is not necessary in theory, as the preconditioners are fixed linear operators, but it is a practical choice to obtain consistent results in the presence of round-off errors. The iteration counts are proxies for the condition number of the preconditioned system, as established by the convergence theory of Krylov methods. 

\begin{table}[htb]
\centering
\setlength{\tabcolsep}{2.5pt}
\begingroup\scriptsize
\setlength{\aboverulesep}{0.2ex}\setlength{\belowrulesep}{0.3ex}
{\small\textbf{(a)} Number of cells}\par\smallskip
\begin{tabular}{ll!{\color{lightgray}\vrule width 0.4pt}rrrrrr}
\toprule
\multicolumn{2}{r}{$\ell$} & 0 & 1 & 2 & 3 & 4 & 5 \\
\midrule
\multirow{3}{*}{2D Cartesian} & Coarse & 16\,384 & 4\,096 & 1\,024 & 256 & 64 & --- \\
 & Medium & 65\,536 & 16\,384 & 4\,096 & 1\,024 & 256 & --- \\
 & Fine & 262\,144 & 65\,536 & 16\,384 & 4\,096 & 1\,024 & 256 \\
\midrule
\multirow{2}{*}{2D simplexified} & Coarse & 20\,000 & 5\,041 & 1\,296 & 324 & 81 & 25 \\
 & Fine & 80\,000 & 20\,164 & 5\,041 & 1\,296 & 324 & 81 \\
\midrule
\multirow{3}{*}{Australia} & Coarse & 8\,981 & 585 & 50 & --- & --- & --- \\
 & Medium & 47\,668 & 3\,138 & 231 & 16 & --- & --- \\
 & Fine & 454\,035 & 28\,773 & 1\,861 & 121 & 9 & --- \\
\midrule
\multirow{2}{*}{3D Cartesian} & Coarse & 4\,096 & 512 & 64 & --- & --- & --- \\
 & Fine & 32\,768 & 4\,096 & 512 & 64 & --- & --- \\
\midrule
\multirow{3}{*}{3D simplexified} & Coarse & 1\,296 & 216 & 27 & --- & --- & --- \\
 & Medium & 4\,374 & 729 & 125 & --- & --- & --- \\
 & Fine & 7\,986 & 1\,000 & 125 & --- & --- & --- \\
\bottomrule
\end{tabular}
\par\medskip
{\small\textbf{(b)} Number of DoFs}\par\smallskip
\begin{tabular}{lll!{\color{lightgray}\vrule width 0.4pt}r!{\color{lightgray}\vrule width 0.4pt}rrrrr!{\color{lightgray}\vrule width 0.4pt}rrrrr}
\toprule
&& & \multicolumn{1}{c}{CG} & \multicolumn{5}{c}{DG} & \multicolumn{5}{c}{HHO} \\
\cmidrule(lr){4-4}\cmidrule(lr){5-9}\cmidrule(lr){10-14}
\multicolumn{3}{r}{$\ell$} & 0 & 1 & 2 & 3 & 4 & 5 & 1 & 2 & 3 & 4 & 5 \\
\midrule
\multirow{6}{*}{2D Cartesian} & \multirow{2}{*}{Coarse} & $\k = 1$ & 16\,129 & 12\,288 & 3\,072 & 768 & 192 & --- & 16\,128 & 3\,968 & 960 & 224 & --- \\
 &  & $\k = 3$ & 146\,689 & 40\,960 & 10\,240 & 2\,560 & 640 & --- & 32\,256 & 7\,936 & 1\,920 & 448 & --- \\
 & \multirow{2}{*}{Medium} & $\k = 1$ & 65\,025 & 49\,152 & 12\,288 & 3\,072 & 768 & --- & 65\,024 & 16\,128 & 3\,968 & 960 & --- \\
 &  & $\k = 3$ & 588\,289 & 163\,840 & 40\,960 & 10\,240 & 2\,560 & --- & 130\,048 & 32\,256 & 7\,936 & 1\,920 & --- \\
 & \multirow{2}{*}{Fine} & $\k = 1$ & 261\,121 & 196\,608 & 49\,152 & 12\,288 & 3\,072 & 768 & 261\,120 & 65\,024 & 16\,128 & 3\,968 & 960 \\
 &  & $\k = 3$ & 2\,356\,225 & 655\,360 & 163\,840 & 40\,960 & 10\,240 & 2\,560 & 522\,240 & 130\,048 & 32\,256 & 7\,936 & 1\,920 \\
\midrule
\multirow{4}{*}{2D simplexified} & \multirow{2}{*}{Coarse} & $\k = 1$ & 9\,801 & 15\,123 & 3\,888 & 972 & 243 & 75 & 22\,358 & 7\,190 & 2\,444 & 779 & 249 \\
 &  & $\k = 3$ & 89\,401 & 50\,410 & 12\,960 & 3\,240 & 810 & 250 & 44\,716 & 14\,508 & 5\,154 & 1\,813 & 669 \\
 & \multirow{2}{*}{Fine} & $\k = 1$ & 39\,601 & 60\,492 & 15\,123 & 3\,888 & 972 & 243 & 84\,610 & 24\,483 & 8\,735 & 3\,058 & 960 \\
 &  & $\k = 3$ & 358\,801 & 201\,640 & 50\,410 & 12\,960 & 3\,240 & 810 & 169\,220 & 49\,055 & 17\,919 & 6\,778 & 2\,426 \\
\midrule
\multirow{6}{*}{Australia} & \multirow{2}{*}{Coarse} & $\k = 1$ & 4\,328 & 1\,755 & 150 & --- & --- & --- & 6\,326 & 674 & --- & --- & --- \\
 &  & $\k = 3$ & 39\,925 & 5\,850 & 500 & --- & --- & --- & 13\,259 & 2\,070 & --- & --- & --- \\
 & \multirow{2}{*}{Medium} & $\k = 1$ & 23\,397 & 9\,414 & 693 & 48 & --- & --- & 34\,871 & 3\,554 & 192 & --- & --- \\
 &  & $\k = 3$ & 213\,193 & 31\,380 & 2\,310 & 160 & --- & --- & 72\,671 & 10\,985 & 633 & --- & --- \\
 & \multirow{2}{*}{Fine} & $\k = 1$ & 225\,415 & 86\,319 & 5\,583 & 363 & 27 & --- & 326\,083 & 30\,925 & 1\,873 & 96 & --- \\
 &  & $\k = 3$ & 2\,038\,348 & 287\,730 & 18\,610 & 1\,210 & 90 & --- & 667\,219 & 94\,570 & 6\,112 & 312 & --- \\
\midrule
\multirow{4}{*}{3D Cartesian} & \multirow{2}{*}{Coarse} & $\k = 1$ & 3\,375 & 2\,048 & 256 & --- & --- & --- & 4\,032 & 432 & --- & --- & --- \\
 &  & $\k = 2$ & 29\,791 & 5\,120 & 640 & --- & --- & --- & 8\,064 & 864 & --- & --- & --- \\
 & \multirow{2}{*}{Fine} & $\k = 1$ & 29\,791 & 16\,384 & 2\,048 & 256 & --- & --- & 34\,560 & 4\,032 & 432 & --- & --- \\
 &  & $\k = 2$ & 250\,047 & 40\,960 & 5\,120 & 640 & --- & --- & 69\,120 & 8\,064 & 864 & --- & --- \\
\midrule
\multirow{6}{*}{3D simplexified} & \multirow{2}{*}{Coarse} & $\k = 1$ & 125 & 864 & 108 & --- & --- & --- & 1\,620 & 273 & --- & --- & --- \\
 &  & $\k = 2$ & 1\,331 & 2\,160 & 270 & --- & --- & --- & 3\,240 & 549 & --- & --- & --- \\
 & \multirow{2}{*}{Medium} & $\k = 1$ & 512 & 2\,916 & 500 & --- & --- & --- & 5\,832 & 900 & --- & --- & --- \\
 &  & $\k = 2$ & 4\,913 & 7\,290 & 1\,250 & --- & --- & --- & 11\,664 & 1\,800 & --- & --- & --- \\
 & \multirow{2}{*}{Fine} & $\k = 1$ & 1\,000 & 4\,000 & 500 & --- & --- & --- & 8\,286 & 1\,098 & --- & --- & --- \\
 &  & $\k = 2$ & 9\,261 & 10\,000 & 1\,250 & --- & --- & --- & 16\,572 & 2\,196 & --- & --- & --- \\
\bottomrule
\end{tabular}
\endgroup
\caption{Mesh hierarchies used throughout \sect{sec:results}. Level $\ell = 0$ is the
fine mesh and $\ell \geq 1$ the agglomerated levels, so the auxiliary space of
the two-level method sits at the $\ell$ whose coarsening ratio is quoted in
Table \ref{tab:niter-two-level}, and the coarsest level of an
$n_{\mathrm{lev}}$-level hierarchy of Table \ref{tab:niter-multilevel} is
$\ell = n_{\mathrm{lev}} - 1$. CG DoFs are those of the conforming fine space,
DG and HHO DoFs those of the coarse space built on each agglomerated level.
Entries marked ``---'' are levels the hierarchy does not provide.}
\label{tab:meshes}
\end{table}

\subsection{Two-level auxiliary space preconditioners}
\label{sec:results-two-level}

For each mesh hierarchy, we consider the finest mesh $\Th$ and take $\TH$ to be any of the coarser meshes $\mathcal{T}_\ell$, $\ell \geq 1$, in the hierarchy, so that the coarsening ratio is the product of the coarsening ratios of the intermediate levels. By doing this, we create a series of two-level problems with increasing coarsening ratios.

We test both the additive preconditioner $P_{\mathrm{add}}$ of \eqref{eq:auxspace-preconditioner} and a multiplicative variant $P_{\mathrm{mult}}$, with the point-Jacobi and vertex-star patch (V-Patch) smoothers described in Section \ref{sec:smoother}. In both cases the smoothing step consists of $n_s$ iterations of the associated damped Richardson smoother $S_{h,n_s,\omega}$.
For $P_{\mathrm{add}}$ we take $n_s = 1$ and $\omega = 1$, which recovers $\Dh^{-1}$ and matches the analysis of Section \ref{sec:smoother}. For $P_{\mathrm{mult}}$ we take $n_s = 5$ and $\omega = 0.1$ in two dimensions and $\omega = 0.05$ in three dimensions, as is common practice for \ac{gmg} methods. In order to isolate the two-level behaviour of the preconditioner, we use a sparse LU factorisation to solve the coarse polytopal problem exactly, that is we set $\BH = \AH^{-1}$.

We report the iteration counts for all five mesh families and both coarse discretisations in Table \ref{tab:niter-two-level}. The columns are labelled by the nominal $\H/\h = \theta^{-\ell}$. At $\k = 1$ the vertex-star patches reduce to single-\ac{dof} blocks, so the two smoothers are the same operator and a single smoother block is reported.

\begin{table}[htb]
\centering
\setlength{\tabcolsep}{2pt}
\begingroup\scriptsize
\setlength{\aboverulesep}{0.2ex}\setlength{\belowrulesep}{0.3ex}
{\small\textbf{(a)} 2D Cartesian meshes}\par\smallskip
\begin{tabular}{ll!{\color{lightgray}\vrule width 0.4pt}ccccc!{\color{lightgray}\vrule width 0.4pt}ccccc!{\color{lightgray}\vrule width 0.4pt}ccccc!{\color{lightgray}\vrule width 0.4pt}ccccc!{\color{lightgray}\vrule width 0.4pt}ccccc!{\color{lightgray}\vrule width 0.4pt}ccccc}
\toprule
& & \multicolumn{15}{c}{CG-DG} & \multicolumn{15}{c}{CG-HHO} \\
\cmidrule(lr){3-17}\cmidrule(lr){18-32}
& & \multicolumn{5}{c}{\shortstack{$\k = 1$\\Jacobi}} & \multicolumn{5}{c}{\shortstack{$\k = 3$\\Jacobi}} & \multicolumn{5}{c}{\shortstack{$\k = 3$\\V-Patch}} & \multicolumn{5}{c}{\shortstack{$\k = 1$\\Jacobi}} & \multicolumn{5}{c}{\shortstack{$\k = 3$\\Jacobi}} & \multicolumn{5}{c}{\shortstack{$\k = 3$\\V-Patch}} \\
\cmidrule(lr){3-7}\cmidrule(lr){8-12}\cmidrule(lr){13-17}\cmidrule(lr){18-22}\cmidrule(lr){23-27}\cmidrule(lr){28-32}
\multicolumn{2}{r}{$\H/\h$} & 2 & 4 & 8 & 16 & 32 & 2 & 4 & 8 & 16 & 32 & 2 & 4 & 8 & 16 & 32 & 2 & 4 & 8 & 16 & 32 & 2 & 4 & 8 & 16 & 32 & 2 & 4 & 8 & 16 & 32 \\
\midrule
\multirow{3}{*}{$P_{\mathrm{add}}$} & Coarse & 12 & 19 & 32 & 53 & --- & 16 & 31 & 67 & 143 & --- & 13 & 16 & 27 & 51 & --- & 3 & 8 & 18 & 48 & --- & 18 & 41 & 91 & 225 & --- & 13 & 18 & 29 & 65 & --- \\
 & Medium & 11 & 19 & 31 & 54 & --- & 14 & 29 & 61 & 129 & --- & 13 & 15 & 25 & 47 & --- & 3 & 7 & 17 & 47 & --- & 18 & 40 & 89 & 237 & --- & 13 & 18 & 29 & 63 & --- \\
 & Fine & 11 & 18 & 30 & 53 & 98 & 13 & 27 & 54 & 113 & 243 & 12 & 14 & 23 & 42 & 87 & 3 & 7 & 16 & 45 & 132 & 18 & 39 & 86 & 235 & 632 & 12 & 18 & 28 & 62 & 161 \\
\midrule
\multirow{3}{*}{$P_{\mathrm{mult}}$} & Coarse & 10 & 16 & 24 & 40 & --- & 11 & 22 & 45 & 95 & --- & 7 & 8 & 14 & 27 & --- & 3 & 7 & 14 & 37 & --- & 14 & 28 & 73 & 190 & --- & 8 & 10 & 16 & 30 & --- \\
 & Medium & 9 & 15 & 25 & 41 & --- & 10 & 20 & 42 & 87 & --- & 7 & 8 & 13 & 24 & --- & 3 & 6 & 14 & 37 & --- & 13 & 28 & 72 & 200 & --- & 8 & 10 & 16 & 29 & --- \\
 & Fine & 8 & 14 & 24 & 42 & 77 & 9 & 18 & 37 & 77 & 167 & 6 & 8 & 12 & 21 & 42 & 3 & 6 & 12 & 37 & 107 & 13 & 28 & 70 & 197 & 549 & 7 & 10 & 16 & 29 & 76 \\
\bottomrule
\end{tabular}
\par\medskip
{\small\textbf{(b)} 2D simplexified meshes}\par\smallskip
\begin{tabular}{ll!{\color{lightgray}\vrule width 0.4pt}ccccc!{\color{lightgray}\vrule width 0.4pt}ccccc!{\color{lightgray}\vrule width 0.4pt}ccccc!{\color{lightgray}\vrule width 0.4pt}ccccc!{\color{lightgray}\vrule width 0.4pt}ccccc!{\color{lightgray}\vrule width 0.4pt}ccccc}
\toprule
& & \multicolumn{15}{c}{CG-DG} & \multicolumn{15}{c}{CG-HHO} \\
\cmidrule(lr){3-17}\cmidrule(lr){18-32}
& & \multicolumn{5}{c}{\shortstack{$\k = 1$\\Jacobi}} & \multicolumn{5}{c}{\shortstack{$\k = 3$\\Jacobi}} & \multicolumn{5}{c}{\shortstack{$\k = 3$\\V-Patch}} & \multicolumn{5}{c}{\shortstack{$\k = 1$\\Jacobi}} & \multicolumn{5}{c}{\shortstack{$\k = 3$\\Jacobi}} & \multicolumn{5}{c}{\shortstack{$\k = 3$\\V-Patch}} \\
\cmidrule(lr){3-7}\cmidrule(lr){8-12}\cmidrule(lr){13-17}\cmidrule(lr){18-22}\cmidrule(lr){23-27}\cmidrule(lr){28-32}
\multicolumn{2}{r}{$\H/\h$} & 2 & 4 & 8 & 16 & 32 & 2 & 4 & 8 & 16 & 32 & 2 & 4 & 8 & 16 & 32 & 2 & 4 & 8 & 16 & 32 & 2 & 4 & 8 & 16 & 32 & 2 & 4 & 8 & 16 & 32 \\
\midrule
\multirow{2}{*}{$P_{\mathrm{add}}$} & Coarse & 20 & 30 & 44 & 76 & 123 & 28 & 45 & 81 & 162 & 283 & 18 & 22 & 32 & 54 & 89 & 17 & 28 & 60 & 110 & 140 & 49 & 121 & 284 & 425 & 420 & 21 & 44 & 100 & 181 & 219 \\
 & Fine & 20 & 31 & 46 & 78 & 151 & 28 & 47 & 82 & 156 & 334 & 18 & 23 & 34 & 53 & 108 & 16 & 29 & 65 & 145 & 262 & 42 & 127 & 326 & 690 & 719 & 20 & 44 & 110 & 251 & 379 \\
\midrule
\multirow{2}{*}{$P_{\mathrm{mult}}$} & Coarse & 18 & 28 & 41 & 56 & 84 & 23 & 41 & 74 & 148 & 206 & 14 & 16 & 25 & 44 & 70 & 21 & 30 & 56 & 98 & 127 & 47 & 104 & 247 & 386 & 341 & 21 & 28 & 56 & 111 & 162 \\
 & Fine & 18 & 29 & 45 & 68 & 110 & 23 & 43 & 76 & 143 & 300 & 14 & 17 & 25 & 43 & 84 & 22 & 31 & 61 & 128 & 237 & 41 & 109 & 282 & 612 & 672 & 21 & 28 & 61 & 149 & 251 \\
\bottomrule
\end{tabular}
\par\medskip
{\small\textbf{(c)} Australia meshes}\par\smallskip
\begin{tabular}{ll!{\color{lightgray}\vrule width 0.4pt}cccc!{\color{lightgray}\vrule width 0.4pt}cccc!{\color{lightgray}\vrule width 0.4pt}cccc!{\color{lightgray}\vrule width 0.4pt}cccc!{\color{lightgray}\vrule width 0.4pt}cccc!{\color{lightgray}\vrule width 0.4pt}cccc}
\toprule
& & \multicolumn{12}{c}{CG-DG} & \multicolumn{12}{c}{CG-HHO} \\
\cmidrule(lr){3-14}\cmidrule(lr){15-26}
& & \multicolumn{4}{c}{\shortstack{$\k = 1$\\Jacobi}} & \multicolumn{4}{c}{\shortstack{$\k = 3$\\Jacobi}} & \multicolumn{4}{c}{\shortstack{$\k = 3$\\V-Patch}} & \multicolumn{4}{c}{\shortstack{$\k = 1$\\Jacobi}} & \multicolumn{4}{c}{\shortstack{$\k = 3$\\Jacobi}} & \multicolumn{4}{c}{\shortstack{$\k = 3$\\V-Patch}} \\
\cmidrule(lr){3-6}\cmidrule(lr){7-10}\cmidrule(lr){11-14}\cmidrule(lr){15-18}\cmidrule(lr){19-22}\cmidrule(lr){23-26}
\multicolumn{2}{r}{$\H/\h$} & 4 & 16 & 64 & 256 & 4 & 16 & 64 & 256 & 4 & 16 & 64 & 256 & 4 & 16 & 64 & 256 & 4 & 16 & 64 & 256 & 4 & 16 & 64 & 256 \\
\midrule
\multirow{3}{*}{$P_{\mathrm{add}}$} & Coarse & 22 & 45 & --- & --- & 44 & 87 & --- & --- & 21 & 34 & --- & --- & 21 & 67 & --- & --- & 87 & 222 & --- & --- & 36 & 105 & --- & --- \\
 & Medium & 21 & 44 & 117 & --- & 59 & 114 & 257 & --- & 25 & 41 & 94 & --- & 23 & 77 & 136 & --- & 113 & 385 & 333 & --- & 42 & 146 & 196 & --- \\
 & Fine & 27 & 61 & 175 & 379 & 61 & 176 & 480 & 1141 & 29 & 62 & 160 & 383 & 26 & 123 & 447 & 336 & 153 & 702 & 1159 & 1027 & 56 & 260 & 641 & 370 \\
\midrule
\multirow{3}{*}{$P_{\mathrm{mult}}$} & Coarse & 22 & 35 & --- & --- & 36 & 81 & --- & --- & 14 & 26 & --- & --- & 23 & 57 & --- & --- & 67 & 190 & --- & --- & 19 & 57 & --- & --- \\
 & Medium & 21 & 40 & 92 & --- & 46 & 105 & 202 & --- & 16 & 30 & 57 & --- & 24 & 65 & 116 & --- & 88 & 320 & 265 & --- & 22 & 75 & 146 & --- \\
 & Fine & 27 & 60 & 147 & 278 & 47 & 162 & 406 & 763 & 18 & 44 & 106 & 242 & 27 & 103 & 392 & 251 & 119 & 573 & 1076 & 687 & 26 & 132 & 397 & 238 \\
\bottomrule
\end{tabular}
\par\medskip
{\small\textbf{(d)} 3D Cartesian meshes}\par\smallskip
\begin{tabular}{ll!{\color{lightgray}\vrule width 0.4pt}ccc!{\color{lightgray}\vrule width 0.4pt}ccc!{\color{lightgray}\vrule width 0.4pt}ccc!{\color{lightgray}\vrule width 0.4pt}ccc!{\color{lightgray}\vrule width 0.4pt}ccc!{\color{lightgray}\vrule width 0.4pt}ccc}
\toprule
& & \multicolumn{9}{c}{CG-DG} & \multicolumn{9}{c}{CG-HHO} \\
\cmidrule(lr){3-11}\cmidrule(lr){12-20}
& & \multicolumn{3}{c}{\shortstack{$\k = 1$\\Jacobi}} & \multicolumn{3}{c}{\shortstack{$\k = 2$\\Jacobi}} & \multicolumn{3}{c}{\shortstack{$\k = 2$\\V-Patch}} & \multicolumn{3}{c}{\shortstack{$\k = 1$\\Jacobi}} & \multicolumn{3}{c}{\shortstack{$\k = 2$\\Jacobi}} & \multicolumn{3}{c}{\shortstack{$\k = 2$\\V-Patch}} \\
\cmidrule(lr){3-5}\cmidrule(lr){6-8}\cmidrule(lr){9-11}\cmidrule(lr){12-14}\cmidrule(lr){15-17}\cmidrule(lr){18-20}
\multicolumn{2}{r}{$\H/\h$} & 2 & 4 & 8 & 2 & 4 & 8 & 2 & 4 & 8 & 2 & 4 & 8 & 2 & 4 & 8 & 2 & 4 & 8 \\
\midrule
\multirow{2}{*}{$P_{\mathrm{add}}$} & Coarse & 13 & 17 & --- & 18 & 28 & --- & 28 & 33 & --- & 14 & 14 & --- & 17 & 31 & --- & 28 & 31 & --- \\
 & Fine & 13 & 19 & 29 & 18 & 27 & 50 & 28 & 32 & 48 & 14 & 16 & 27 & 17 & 31 & 67 & 28 & 32 & 49 \\
\midrule
\multirow{2}{*}{$P_{\mathrm{mult}}$} & Coarse & 11 & 13 & --- & 18 & 25 & --- & 15 & 17 & --- & 14 & 14 & --- & 18 & 31 & --- & 15 & 18 & --- \\
 & Fine & 13 & 15 & 23 & 19 & 29 & 49 & 15 & 17 & 25 & 13 & 16 & 29 & 18 & 34 & 76 & 15 & 17 & 30 \\
\bottomrule
\end{tabular}
\par\medskip
{\small\textbf{(e)} 3D simplexified meshes}\par\smallskip
\begin{tabular}{ll!{\color{lightgray}\vrule width 0.4pt}cc!{\color{lightgray}\vrule width 0.4pt}cc!{\color{lightgray}\vrule width 0.4pt}cc!{\color{lightgray}\vrule width 0.4pt}cc!{\color{lightgray}\vrule width 0.4pt}cc!{\color{lightgray}\vrule width 0.4pt}cc}
\toprule
& & \multicolumn{6}{c}{CG-DG} & \multicolumn{6}{c}{CG-HHO} \\
\cmidrule(lr){3-8}\cmidrule(lr){9-14}
& & \multicolumn{2}{c}{\shortstack{$\k = 1$\\Jacobi}} & \multicolumn{2}{c}{\shortstack{$\k = 2$\\Jacobi}} & \multicolumn{2}{c}{\shortstack{$\k = 2$\\V-Patch}} & \multicolumn{2}{c}{\shortstack{$\k = 1$\\Jacobi}} & \multicolumn{2}{c}{\shortstack{$\k = 2$\\Jacobi}} & \multicolumn{2}{c}{\shortstack{$\k = 2$\\V-Patch}} \\
\cmidrule(lr){3-4}\cmidrule(lr){5-6}\cmidrule(lr){7-8}\cmidrule(lr){9-10}\cmidrule(lr){11-12}\cmidrule(lr){13-14}
\multicolumn{2}{r}{$\H/\h$} & 2 & 4 & 2 & 4 & 2 & 4 & 2 & 4 & 2 & 4 & 2 & 4 \\
\midrule
\multirow{3}{*}{$P_{\mathrm{add}}$} & Coarse & 10 & 17 & 24 & 29 & 18 & 20 & 19 & 30 & 48 & 67 & 39 & 46 \\
 & Medium & 13 & 16 & 26 & 33 & 20 & 22 & 29 & 35 & 54 & 90 & 43 & 57 \\
 & Fine & 22 & 21 & 39 & 43 & 26 & 27 & 54 & 60 & 84 & 123 & 59 & 79 \\
\midrule
\multirow{3}{*}{$P_{\mathrm{mult}}$} & Coarse & 11 & 14 & 28 & 29 & 22 & 18 & 24 & 33 & 64 & 77 & 50 & 52 \\
 & Medium & 14 & 14 & 35 & 38 & 27 & 23 & 37 & 44 & 74 & 112 & 58 & 71 \\
 & Fine & 26 & 19 & 52 & 51 & 35 & 33 & 71 & 75 & 111 & 148 & 79 & 96 \\
\bottomrule
\end{tabular}
\endgroup
\caption{Two-level preconditioners: iteration counts for the mesh families of
Table \ref{tab:meshes}, by coarse space, polynomial degree, smoother and
nominal coarsening ratio $\H/\h$. ``---'' marks ratios the hierarchy does
not provide.}
\label{tab:niter-two-level}
\end{table}

On the Cartesian families, where mesh regularity is preserved exactly, the results match the theory. For a fixed coarsening ratio the iteration counts are independent of the fine mesh, and they grow mildly and monotonically with $\H/\h$, well below the growth allowed by the bounds of Sections \ref{sec:dg} and \ref{sec:hho}. On the unstructured families the behaviour is qualitatively the same, but the counts depend on the quality of the agglomeration. In particular, on the \emph{Fine} Australian hierarchy the agglomeration quality at $\ell = 3$ is noticeably poorer, and as a consequence the iteration counts for \ac{hho} are significantly higher than at the other levels, including $\ell = 4$. Comparing the smoothers, the iteration counts with the vertex-star patch smoother are insensitive to the polynomial degree, whereas point-Jacobi degrades markedly as the degree grows. Such $p$-robustness of patch smoothers is well documented for conforming discretisations \cite{10.1093/imanum/drl046,Pavarino1993}, but has not, to our knowledge, been analysed for our particular setting. The multiplicative variant $P_{\mathrm{mult}}$ generally, but not uniformly, needs fewer iterations than $P_{\mathrm{add}}$.

\subsection{Multilevel coarse solvers}
\label{sec:results-multilevel}

We now replace the exact coarse solve by a single iteration of the \ac{gmg} V-cycle described in Section \ref{sec:multilevel}.
On the coarser \ac{gmg} levels, we use the damped Richardson smoother $S_{\ell,n_s,\omega}$ with $n_s = 5$ and $\omega = 0.1$ in two dimensions and $\omega = 0.05$ in three. 

\begin{table}[htb]
\centering
\setlength{\tabcolsep}{2pt}
\begingroup\scriptsize
\setlength{\aboverulesep}{0.2ex}\setlength{\belowrulesep}{0.3ex}
\begin{tabular}[c]{@{}c@{}}
{\small\textbf{(a)} 2D Cartesian meshes, $\H/\h = 2$}\\[\smallskipamount]
\begin{tabular}{ll!{\color{lightgray}\vrule width 0.4pt}ccccc!{\color{lightgray}\vrule width 0.4pt}ccccc!{\color{lightgray}\vrule width 0.4pt}ccccc!{\color{lightgray}\vrule width 0.4pt}ccccc}
\toprule
& & \multicolumn{10}{c}{CG-DG} & \multicolumn{10}{c}{CG-HHO} \\
\cmidrule(lr){3-12}\cmidrule(lr){13-22}
& & \multicolumn{5}{c}{$\k = 1$} & \multicolumn{5}{c}{$\k = 2$} & \multicolumn{5}{c}{$\k = 1$} & \multicolumn{5}{c}{$\k = 2$} \\
\cmidrule(lr){3-7}\cmidrule(lr){8-12}\cmidrule(lr){13-17}\cmidrule(lr){18-22}
\multicolumn{2}{r}{$n_{\mathrm{lev}}$} & 2 & 3 & 4 & 5 & 6 & 2 & 3 & 4 & 5 & 6 & 2 & 3 & 4 & 5 & 6 & 2 & 3 & 4 & 5 & 6 \\
\midrule
\multirow{3}{*}{$P_{\mathrm{add}}$} & Coarse & 12 & 14 & 16 & --- & --- & 14 & 16 & 18 & --- & --- & 3 & 10 & 12 & --- & --- & 10 & 17 & 19 & --- & --- \\
 & Medium & 11 & 14 & 15 & 16 & --- & 13 & 16 & 17 & 19 & --- & 3 & 9 & 11 & 13 & --- & 10 & 16 & 19 & 20 & --- \\
 & Fine & 11 & 13 & 14 & 16 & 17 & 13 & 14 & 17 & 18 & 19 & 3 & 8 & 10 & 13 & 13 & 9 & 15 & 19 & 19 & 20 \\
\midrule
\multirow{3}{*}{$P_{\mathrm{mult}}$} & Coarse & 10 & 12 & 12 & --- & --- & 7 & 8 & 9 & --- & --- & 3 & 7 & 8 & --- & --- & 7 & 9 & 9 & --- & --- \\
 & Medium & 9 & 11 & 12 & 13 & --- & 7 & 8 & 9 & 9 & --- & 3 & 6 & 8 & 9 & --- & 7 & 9 & 9 & 10 & --- \\
 & Fine & 8 & 10 & 12 & 13 & 14 & 6 & 7 & 9 & 9 & 9 & 3 & 6 & 8 & 9 & 9 & 7 & 9 & 9 & 10 & 10 \\
\bottomrule
\end{tabular}
\end{tabular}
\par\medskip
\begin{tabular}[c]{@{}c@{}}
{\small\textbf{(b)} 2D Cartesian meshes, $\H/\h = 4$}\\[\smallskipamount]
\begin{tabular}{ll!{\color{lightgray}\vrule width 0.4pt}ccc!{\color{lightgray}\vrule width 0.4pt}ccc!{\color{lightgray}\vrule width 0.4pt}ccc!{\color{lightgray}\vrule width 0.4pt}ccc}
\toprule
& & \multicolumn{6}{c}{CG-DG} & \multicolumn{6}{c}{CG-HHO} \\
\cmidrule(lr){3-8}\cmidrule(lr){9-14}
& & \multicolumn{3}{c}{$\k = 1$} & \multicolumn{3}{c}{$\k = 2$} & \multicolumn{3}{c}{$\k = 1$} & \multicolumn{3}{c}{$\k = 2$} \\
\cmidrule(lr){3-5}\cmidrule(lr){6-8}\cmidrule(lr){9-11}\cmidrule(lr){12-14}
\multicolumn{2}{r}{$n_{\mathrm{lev}}$} & 2 & 3 & 4 & 2 & 3 & 4 & 2 & 3 & 4 & 2 & 3 & 4 \\
\midrule
\multirow{3}{*}{$P_{\mathrm{add}}$} & Coarse & 19 & 24 & --- & 19 & 26 & --- & 8 & 19 & --- & 16 & 29 & --- \\
 & Medium & 19 & 25 & 26 & 18 & 25 & 28 & 7 & 18 & 20 & 16 & 28 & 30 \\
 & Fine & 18 & 24 & 27 & 17 & 24 & 27 & 7 & 17 & 20 & 15 & 27 & 30 \\
\midrule
\multirow{3}{*}{$P_{\mathrm{mult}}$} & Coarse & 16 & 18 & --- & 10 & 14 & --- & 7 & 15 & --- & 10 & 15 & --- \\
 & Medium & 15 & 19 & 19 & 10 & 13 & 15 & 6 & 14 & 16 & 10 & 15 & 16 \\
 & Fine & 14 & 19 & 20 & 9 & 13 & 14 & 6 & 14 & 16 & 10 & 14 & 16 \\
\bottomrule
\end{tabular}
\end{tabular}
\hfill
\begin{tabular}[c]{@{}c@{}}
{\small\textbf{(d)} 3D Cartesian meshes, $\H/\h = 2$}\\[\smallskipamount]
\begin{tabular}{ll!{\color{lightgray}\vrule width 0.4pt}ccc!{\color{lightgray}\vrule width 0.4pt}ccc!{\color{lightgray}\vrule width 0.4pt}ccc!{\color{lightgray}\vrule width 0.4pt}ccc}
\toprule
& & \multicolumn{6}{c}{CG-DG} & \multicolumn{6}{c}{CG-HHO} \\
\cmidrule(lr){3-8}\cmidrule(lr){9-14}
& & \multicolumn{3}{c}{$\k = 1$} & \multicolumn{3}{c}{$\k = 2$} & \multicolumn{3}{c}{$\k = 1$} & \multicolumn{3}{c}{$\k = 2$} \\
\cmidrule(lr){3-5}\cmidrule(lr){6-8}\cmidrule(lr){9-11}\cmidrule(lr){12-14}
\multicolumn{2}{r}{$n_{\mathrm{lev}}$} & 2 & 3 & 4 & 2 & 3 & 4 & 2 & 3 & 4 & 2 & 3 & 4 \\
\midrule
\multirow{2}{*}{$P_{\mathrm{add}}$} & Coarse & 13 & 14 & --- & 28 & 29 & --- & 14 & 14 & --- & 28 & 29 & --- \\
 & Fine & 13 & 16 & 16 & 28 & 29 & 30 & 14 & 14 & 15 & 28 & 28 & 29 \\
\midrule
\multirow{2}{*}{$P_{\mathrm{mult}}$} & Coarse & 11 & 10 & --- & 15 & 15 & --- & 14 & 15 & --- & 15 & 15 & --- \\
 & Fine & 13 & 12 & 12 & 15 & 15 & 15 & 13 & 15 & 15 & 15 & 15 & 15 \\
\bottomrule
\end{tabular}
\end{tabular}
\par\medskip
\begin{tabular}[c]{@{}c@{}}
{\small\textbf{(c)} Australia meshes, $\H/\h \simeq 4$}\\[\smallskipamount]
\begin{tabular}{ll!{\color{lightgray}\vrule width 0.4pt}cccc!{\color{lightgray}\vrule width 0.4pt}cccc!{\color{lightgray}\vrule width 0.4pt}cccc!{\color{lightgray}\vrule width 0.4pt}cccc}
\toprule
& & \multicolumn{8}{c}{CG-DG} & \multicolumn{8}{c}{CG-HHO} \\
\cmidrule(lr){3-10}\cmidrule(lr){11-18}
& & \multicolumn{4}{c}{$\k = 1$} & \multicolumn{4}{c}{$\k = 2$} & \multicolumn{4}{c}{$\k = 1$} & \multicolumn{4}{c}{$\k = 2$} \\
\cmidrule(lr){3-6}\cmidrule(lr){7-10}\cmidrule(lr){11-14}\cmidrule(lr){15-18}
\multicolumn{2}{r}{$n_{\mathrm{lev}}$} & 2 & 3 & 4 & 5 & 2 & 3 & 4 & 5 & 2 & 3 & 4 & 5 & 2 & 3 & 4 & 5 \\
\midrule
\multirow{3}{*}{$P_{\mathrm{add}}$} & Coarse & 22 & 26 & --- & --- & 21 & 25 & --- & --- & 21 & 25 & --- & --- & 27 & 27 & --- & --- \\
 & Medium & 21 & 24 & 25 & --- & 24 & 27 & 46 & --- & 23 & 28 & 31 & --- & 32 & 32 & 32 & --- \\
 & Fine & 27 & 30 & 51 & 51 & 28 & 34 & 130 & 170 & 26 & 33 & 36 & 37 & 42 & 41 & 41 & 42 \\
\midrule
\multirow{3}{*}{$P_{\mathrm{mult}}$} & Coarse & 22 & 20 & --- & --- & 17 & 19 & --- & --- & 23 & 27 & --- & --- & 21 & 23 & --- & --- \\
 & Medium & 21 & 22 & 18 & --- & 18 & 19 & 25 & --- & 24 & 29 & 32 & --- & 22 & 25 & 27 & --- \\
 & Fine & 27 & 29 & 41 & 38 & 21 & 24 & 93 & 103 & 27 & 34 & 38 & 40 & 28 & 31 & 34 & 35 \\
\bottomrule
\end{tabular}
\end{tabular}
\hfill
\begin{tabular}[c]{@{}c@{}}
{\small\textbf{(e)} 3D simplexified meshes, $\H/\h \simeq 2$}\\[\smallskipamount]
\begin{tabular}{ll!{\color{lightgray}\vrule width 0.4pt}cc!{\color{lightgray}\vrule width 0.4pt}cc!{\color{lightgray}\vrule width 0.4pt}cc!{\color{lightgray}\vrule width 0.4pt}cc}
\toprule
& & \multicolumn{4}{c}{CG-DG} & \multicolumn{4}{c}{CG-HHO} \\
\cmidrule(lr){3-6}\cmidrule(lr){7-10}
& & \multicolumn{2}{c}{$\k = 1$} & \multicolumn{2}{c}{$\k = 2$} & \multicolumn{2}{c}{$\k = 1$} & \multicolumn{2}{c}{$\k = 2$} \\
\cmidrule(lr){3-4}\cmidrule(lr){5-6}\cmidrule(lr){7-8}\cmidrule(lr){9-10}
\multicolumn{2}{r}{$n_{\mathrm{lev}}$} & 2 & 3 & 2 & 3 & 2 & 3 & 2 & 3 \\
\midrule
\multirow{3}{*}{$P_{\mathrm{add}}$} & Coarse & 10 & 14 & 18 & 15 & 19 & 26 & 39 & 41 \\
 & Medium & 13 & 14 & 20 & 18 & 29 & 32 & 43 & 47 \\
 & Fine & 22 & 19 & 26 & 22 & 54 & 56 & 59 & 62 \\
\midrule
\multirow{3}{*}{$P_{\mathrm{mult}}$} & Coarse & 11 & 10 & 22 & 15 & 24 & 32 & 50 & 53 \\
 & Medium & 14 & 11 & 27 & 21 & 37 & 41 & 58 & 63 \\
 & Fine & 26 & 18 & 35 & 27 & 71 & 71 & 79 & 81 \\
\bottomrule
\end{tabular}
\end{tabular}
\endgroup
\caption{Multilevel preconditioners: iteration counts with the auxiliary space one
agglomeration level above the fine mesh and a \ac{gmg} coarse solver on the
remaining levels, $n_{\mathrm{lev}}$ being the total number of levels, so
that $n_{\mathrm{lev}} = 2$ is the two-level method with an exact coarse
solve. The fine-level smoother is V-Patch throughout. ``---'' marks depths
the hierarchy does not provide.}
\label{tab:niter-multilevel}
\end{table}

Iteration counts for the multilevel preconditioners are reported in Table \ref{tab:niter-multilevel}. Replacing the exact coarse solve by a single \ac{gmg} V-cycle costs a handful of additional iterations. On the Cartesian families in two and three dimensions, and on the three-dimensional simplexified meshes, adding agglomeration levels leaves the counts unchanged, so neither the fine mesh size nor the number of levels affects them and the preconditioners are scalable. As in the two-level case, the quality of the agglomeration has a significant impact on the unstructured families. On the \emph{Fine} Australian hierarchy, the \ac{dg} iteration counts are constant for the two- and three-level hierarchies, increase noticeably when the level $\ell = 3$ is added, and stabilise afterwards, whereas the \ac{hho} variant is unaffected. This confirms that the quality of the agglomeration, rather than the number of levels, is what matters on the unstructured families.

\begin{remark}[On the choice of constants] \label{rem:constants}
  The \ac{dg} penalty parameter $\sigma$ and the smoother weights $\omega$ used throughout this section are, by choice, not optimal. On most experiments the penalty can be taken smaller and the damping weights larger, both improving performance of the solver. To keep the results directly comparable across experiments, we have instead opted for a fixed set of parameters. Tuning these parameters for a given problem would yield noticeably lower iteration counts than the ones reported.
\end{remark}

\section{Conclusions}
\label{sec:conclusions}

We have introduced and analysed a family of multilevel auxiliary space preconditioners for conforming Lagrangian \ac{fe} discretisations of elliptic problems. At the first level, the auxiliary problem is a non-conforming polytopal discretisation posed on an arbitrarily agglomerated coarse mesh, and at the subsequent levels that problem is solved by an agglomeration-based \ac{gmg} method, so that the whole hierarchy is generated from the fine mesh by agglomeration alone. We have considered two polytopal discretisations, an SIPG \ac{dg} method and a statically condensed mixed-order \ac{hho} method, the latter extending directly to \ac{hdg}. In both cases, the prolongation is an element-wise interpolation onto the fine mesh followed by nodal averaging, preceded by a local potential reconstruction in the \ac{hho} case. Under standard regularity assumptions on the agglomerated mesh, we have proven that the preconditioner satisfies $\kappa \lesssim (\H/\h)^3$, and $\kappa \lesssim (\H/\h)^2$ for the \ac{dg} variant with penalty $\sigma \simeq \H/\h$, for point-Jacobi and vertex-star patch smoothers when the coarse problem is solved exactly. The \ac{gmg} coarse solver enters these bounds only through the extreme eigenvalues of the preconditioned coarse problem, which are bounded independently of the number of levels by the results of \cite{badia2026geometricmultigridsolvershybrid,Antonietti2015Multigrid,Antonietti2017Multigrid}, so that the whole solver is scalable.
Numerical experiments on Cartesian, simplexified and unstructured meshes of a complicated geometry, the last two with agglomerates generated by automatic clustering algorithms, corroborate the theory. The quality of the agglomeration has a visible effect on the preconditioner. Well-shaped agglomerates give iteration counts that are independent of the fine mesh size and of the number of levels, whereas poorly shaped ones increase them. Agglomeration strategies that guarantee well-shaped coarse cells are the main open question for the practical use of these preconditioners.

\section*{Acknowledgments}

This research was partially funded by the Australian Government through the Australian Research Council (project numbers DP210103092 and DP220103160). This work was also supported by computational resources provided by the Australian Government through NCI and Pawsey under the NCMAS Merit Allocation Schemes. 

\printbibliography

\end{document}